\documentclass[onefignum,onetabnum,final]{siamart171218}

\newcommand{\pderiv}[2]{\frac{\partial #1}{\partial #2}} 
\newcommand{\red}[1]{{\color{red} #1}}

\newcommand{\bigoh}[1]{\mathcal{O}\left( #1 \right)}
\newcommand\norm[1]{\left\lVert#1\right\rVert}
\newcommand{\chevron}[1]{\left\langle #1 \right\rangle}
\newcommand{\aeps}[0]{a^{\epsilon}}
\newcommand{\Keta}[0]{K_{\eta}}
\newcommand{\rhoeps}[0]{\rho^{\epsilon}}
\newcommand{\intdelta}[0]{\frac{1}{\delta} \int_{-\delta/2}^{\delta/2} \int_{-1/2}^{1/2}}
\DeclareMathOperator{\sign}{sign}

\usepackage{comment}
\usepackage{lipsum}
\usepackage{amssymb}
\usepackage{tikz} 
\usepackage{amsfonts}
\usepackage{mathtools} 
\DeclarePairedDelimiter\ceil{\lceil}{\rceil}
\DeclarePairedDelimiter\floor{\lfloor}{\rfloor}
\usepackage{graphicx}
\usepackage[caption=false]{subfig}
\usepackage{epstopdf}
\usepackage{algorithmic}
\ifpdf
  \DeclareGraphicsExtensions{.eps,.pdf,.png,.jpg}
\else
  \DeclareGraphicsExtensions{.eps}
\fi

\newsiamremark{remark}{Remark}
\newsiamremark{hypothesis}{Hypothesis}
\crefname{hypothesis}{Hypothesis}{Hypotheses}
\newsiamthm{claim}{Claim}

\headers{Reducing the boundary error in numerical homogenization}{Sean P. Carney, Milica Dussinger, and Bj\"{o}rn Engquist}

\title{On the nature of the boundary resonance error in numerical homogenization and its reduction
\thanks{
\funding{This research is supported by National Science Foundation Grant DMS-1620396.}
}}

\author{
Sean P. Carney\footnotemark[1]
\and Milica Dussinger\footnotemark[2] \and
Bj\"{o}rn Engquist\footnotemark[3]
} 

\usepackage{amsopn}

\ifpdf
\hypersetup{
  pdftitle={On the nature of the boundary resonance error in numerical homogenization and its reduction},
  pdfauthor={Sean P. Carney, Milica Dussinger, and Bj\"{o}rn Engquist}
}
\fi

\graphicspath{ {figs/} }
\DeclareGraphicsExtensions{.pdf}
\begin{document}

\maketitle
\renewcommand{\thefootnote}{\fnsymbol{footnote}}
\footnotetext[1]{Department of Mathematical Sciences, George Mason University, Fairfax, VA, 22030, USA; scarney6@gmu.edu}
\footnotetext[2]{milicadussinger@gmail.com}
\footnotetext[3]{Department of Mathematics and Oden Institute for Computational Engineering and Sciences, The University of Texas at Austin, Austin, TX, 78712, USA; engquist@oden.utexas.edu}

\begin{abstract}
Numerical homogenization of multiscale 
equations typically requires taking
an average of the solution to a microscale problem. Both the boundary conditions and domain 
size of the microscale problem play an important role in the accuracy of the homogenization 
procedure. In particular, imposing naive boundary conditions leads to a $\mathcal{O}(\epsilon/\eta)$ 
error in the computation, where $\epsilon$ is the characteristic size of the microscopic 
fluctuations in the heterogeneous media, and $\eta$ is the size of the microscopic domain.
 This so-called boundary, or ``cell resonance" error can dominate discretization error and pollute the entire 
homogenization scheme. 
There exist several techniques in the literature
 to reduce the error. Most strategies involve modifying the form of the microscale cell
problem. Below we present an alternative procedure based on the observation that the resonance
 error itself is an oscillatory function of domain size $\eta$. 
After rigorously characterizing the 
oscillatory behavior for one dimensional and quasi-one dimensional microscale domains, 
we present a novel strategy to 
reduce the resonance error.
Rather than modifying the form of the cell problem, the original problem is solved for a 
sequence of domain sizes, and the results are averaged against kernels satisfying certain moment
conditions and regularity properties.
Numerical examples in one and two dimensions illustrate the utility of the approach. 
%
%

\end{abstract}



\section{Introduction}
\label{sec:intro}
This work is concerned with developing methods for the numerical homogenization of the multiscale elliptic equation
\begin{align}\label{eq:ms_elliptic}
-\nabla \cdot\left( a^{\epsilon} \nabla u^{\epsilon}\right) &= f, \qquad \Omega \subset \mathbb{R}^d\\
u^{\epsilon} &= g, \qquad \partial \Omega \nonumber
\end{align}
where $0 < \epsilon \ll 1$, $|\Omega| = \mathcal{O}(1)$, and $a^{\epsilon}$ is assumed to be symmetric,
 positive-definite, and bounded, so that $\forall w \in \mathbb{R}^d$,  
\begin{equation}\label{eq:a_bounds}
\lambda |w|^2 \le \chevron{w, a^{\epsilon}(x)w} \le \Lambda |w|^2, 
\end{equation}
for some $0 < \lambda < \Lambda < \infty$ and almost every $x \in \Omega$. $a^{\epsilon}$ is also assumed to oscillate rapidly, 
that is, with frequencies $\sim 1/\epsilon$, which in turn implies the solution to the elliptic equation
$u^{\epsilon}$ is oscillatory as well. Equation \eqref{eq:ms_elliptic} models, for instance, 
steady-state heat conduction in a composite material \cite{Zweben:1998}, as well as porous media flow through material 
with highly variable permeability \cite{Durlofsky1998}. The discussion below is also valid for 
numerical homogenization of multiscale parabolic equations, as well as second order hyperbolic equations. 

The mathematical study of the behavior of the solution to (\ref{eq:ms_elliptic}) as $\epsilon \to 0$
has a long history; see \cite{papanicolaou:1978,kozlov:1994,donato:1999}. 
It is well known that when $a^{\epsilon}$ is periodic, $u^{\epsilon} \to \overline{u}$ weakly in 
$H^1(\Omega)$ and strongly in $L^2(\Omega)$, where $\overline{u}$ satisfies the homogenized equation
\begin{align*}\label{eq:hom_elliptic}
-\nabla \cdot \left(\overline{a}\, \nabla\, \overline{u}\right) &= f, \qquad \Omega \\
\overline{u} &= g, \qquad \partial \Omega. \nonumber 
\end{align*} 
The entries in the homogenized tensor $\overline{a}$ are given by 
\begin{equation}\label{eq:hom_tensor}
\overline{a}_{ij} = \frac{1}{|Y|} \int_Y \left(a_{ij}(y) + a_{ik}(y)\frac{\partial \chi_j}{\partial y_k}(y)\right) \, dy,
\end{equation}
where summation over repeated indices is implied, and $\chi_j$ solves the cell problem posed in the unit cell $Y = [-1/2,1/2]^d$
\begin{align}\label{eq:cell_problem}
-\nabla \cdot \left( a \nabla \chi_j \right) &= \nabla \cdot \left(a\, e_j\right), \qquad Y \\
\chi_j \,\,\,\,&Y \text{-periodic and mean zero,} \nonumber
\end{align}
and $e_j$ is the standard unit basis vector in $\mathbb{R}^d$. 
With a slight modification, a similar result holds for locally periodic media. In $d=1$, the cell problem reduces to a 
two-point boundary value problem, and it is easy to show that $\overline{a}$ is given by the harmonic
average of $a$
\begin{equation*}\label{eq:har_avg}
\overline{a} = \left( \int_{-1/2}^{1/2} \frac{1}{a(s)} \, ds \right)^{-1}.
\end{equation*}

In more general settings,
the theory of $\Gamma$-convergence provides sufficient conditions for $u^{\epsilon}$ to limit
to some $\overline{u}$; 
 however, 
there does not in general exist a closed expression for $\overline{a}$, and numerical homogenization techniques
are required to approximate the homogenized tensor.

There exist many strategies for the numerical homogenization of (\ref{eq:ms_elliptic}); examples include 
the multiscale finite element method (MsFEM) \cite{hou:1997,hou:1999,hou:2009}, upscaling techniques \cite{arbogast:2000}, 
the variational multiscale method \cite{hughes:1998}, the heterogeneous 
multiscale method (HMM) \cite{engquist:2007,abdulle:2012}, the equation-free method \cite{kevrekidis2003equation}, and local 
orthogonal decomposition (LOD) \cite{peterseim:2020}. Review articles can be found in 
\cite{gloria_review} and \cite{altmann2021}. 
Each method typically involves solving a microscale problem to estimate some missing data from 
the macroscopic model. Because many examples of realistic heterogeneous media are approximately 
periodic \cite{cartwright2007impact,freudiger2008label}, in practice, both the size of the microscale domain and the 
boundary conditions prescribed for the microscale problem play an important role in the 
accuracy of the homogenization procedure.

A common way to study the effect that the choice of microscale domain size and boundary conditions
has on numerical homogenization techniques is the setting of periodic oscillations, 
where analytic results are known, but where the period is assumed to be unknown. In this context,
E, Ming, and Zhang \cite{pingbingming:2004} first proved that if 
the microscale problem is solved on a domain whose length $\eta$ does not coincide with the true period 
of the heterogeneous media $\epsilon$, a boundary, or ``cell resonance" error results that is proportional to $\epsilon/\eta$ 
and $\epsilon/\eta + \eta$ for truly periodic and locally periodic data, respectively.
Yue and E \cite{yue:2006}
performed a numerical study of this boundary error as a function of domain length $\eta$ for 
Dirichlet, Neumann, and periodic boundary conditions and found results consistent with 
the theory from \cite{pingbingming:2004}.
Since an $\mathcal{O}(\epsilon/\eta)$
error can potentially dominate discretization error in any sensible numerical homogenization scheme, 
its reduction can be of great practical importance. One strategy to reduce the cell resonance error is to 
simply take $\eta \gg \epsilon$. This, however, will of course greatly increase the overall 
computational cost of the homogenization procedure. Instead, it is preferable to devise strategies 
to asymptotically reduce the cell resonance error to $\mathcal{O}(\epsilon/\eta)^r$ for some 
$r > 1$, so that greater accuracy can be obtained at a lower cost.

Several techniques to address the cell resonance error exist in the literature. 
Most involve 
the use of smooth, compactly supported averaging kernels, or ``masks'', when computing the 
homogenized coefficients; their compact support helps reduce the 
influence of the boundary error. The masks are typically inserted into the ``data estimation'', or 
upscaling step of the numerical homogenization procedure, that is, the method-specific analog of \eqref{eq:hom_tensor}.
 In the context of the MsFEM, for example, this is called oversampling 
\cite{peterseim:2013,hou:2004}, while in HMM it is called filtering \cite{yue:2006,engquist:2007}.
Although these kernels indeed help reduce the cell resonance error, they only help lower 
the prefactor in the error's rate of decay; the rate itself remains proportional to 
$\epsilon/\eta$ \cite{yue:2006,Gloria:2008}. 
In \cite{blanc:2010}, however, Blanc and Le Bris proposed to insert the averaging kernels directly 
in the cell problem \eqref{eq:cell_problem} itself, in contrast to simply utilizing them 
as post-processing tools. In general this technique helped improve the decay rate of the cell resonance
error to second order in $\epsilon/\eta$. 

There are a few other notable examples that 
modify the form of the cell problem to reduce the influence
of the boundary. In \cite{gloria:2011}, Gloria modified the microscale 
elliptic problem to include a zeroth order term. Since the Green's
function associated to the modified elliptic operator decays exponentially fast, the influence 
from the boundary is lessened, and a convergence rate of $\mathcal{O}(\epsilon/\eta)^4$ for
the cell resonance error was obtained in the periodic setting. The approach was improved in a 
follow up work \cite{gloria:2016} using Richardson extrapolation. 
In \cite{runborg:2016}, Arjmand and Runborg proposed to solve a time-dependent hyperbolic PDE on the
microscale, based on the idea that the finite speed of propagation of the initial data 
prevents the boundary error from affecting the solution in interior of the domain, provided 
the domain is sufficiently large. Using averaging kernels with special regularity and vanishing
moment properties, they were able to rigorously obtain in the periodic setting
 a convergence rate of $\mathcal{O}(\epsilon/\eta)^q$ 
for arbitrary $q$, depending on the regularity properties of the averaging kernel.
More recently, Abdulle et al.\ \cite{abdulle:2019,abdulle2021parabolic,abdulle2023elliptic} proposed two separate
but related methods to reduce the boundary error. The first approach solves
a parabolic microscale problem over a time interval $T$ that depends on the domain length
$\eta$ and the bounds for the oscillatory tensor $a$ \eqref{eq:a_bounds}. Similar to 
\cite{gloria:2011}, the method is 
based on the exponential decay of the parabolic Green's function. The second approach solves 
an elliptic equation with a modified forcing term of the form
\begin{equation}\label{eq:abdulle_forcing}
f = \sum_{k=1}^N e^{-\lambda_k T} \chevron{\nabla\cdot a, \varphi_k} \varphi_k,
\end{equation}
where $\{\lambda_k, \varphi_k\}_{k=1}^N$ are the first $N$ dominant eigenvalues and eigenfunctions
of the operator $A:= -\nabla \cdot \left( a \nabla\right)$ on $[-\eta/2,\eta/2]^d$ equipped with 
Dirichlet boundary conditions; this $f$ is a spectral truncation of the 
solution operator $e^{-AT}$ for the parabolic problem from semigroup theory. With appropriate 
selections for the parameters $T$ and $N$, the authors prove exponential decay of the boundary error; 
see also \cite{abdulle2020mathicse} for an extension to the case of stochastic coefficients. 

Rather than modifying the form of the cell problem, we present 
below an alternative method for reducing the 
cell resonance error based on the simple observation that, as a 
function of the microscale domain size $\eta$, the error itself 
is oscillatory. This  
was observed in the first numerical study of Yue and E; see e.g.\
Fig.~6 in \cite{yue:2006}. 
One goal of the present work is to more systematically characterize
the resonance error as $\eta$ changes. 
In dimension $d=1$ we show that the error can be written 
\begin{equation}\label{eq:cell_res_with_corrector}
\frac{1}{x} P(x) + R(x),
\end{equation}
where $x = \eta/\epsilon$, $P(x)$ is a periodic function, 
and $R(x)$ is a ``corrector'' term given by 
\begin{equation}\label{eq:one_d_corrector}
R(x) = \sum_{k=2}^{\infty} \left(-\frac{1}{x} P(x)\right)^k .
\end{equation}
In the case of dimension $d > 1$, a decomposition of the form 
\eqref{eq:cell_res_with_corrector} is more 
difficult to show, as the resonance error depends on the solution to an elliptic 
problem with boundary data that oscillates as a function of $x$.
The oscillatory boundary 
condition can lead to complicated boundary layers in the interior; see
e.g.\ \cite{masmoudi:2012} and \cite{feldman:2014}. For two dimensional 
microscale domains of the form 
\begin{equation} \label{eq:I_eta_tubular}
I_{\eta} = [-\eta/2, \eta/2]\times [-\epsilon/2,\epsilon/2],
\end{equation}
we show using Floquet theory \cite{floquet1883equations} that the cell resonance 
error can also be written in the form \eqref{eq:cell_res_with_corrector}, where 
the corrector $R(x)$ is given by the average of a locally 
$\epsilon$-periodic function. 
More generally for domains of the form $[-\eta/2,\eta/2]^d$ that are relevant to 
numerical homogenization we offer numerical 
evidence that something analogous to \eqref{eq:cell_res_with_corrector} 
indeed holds, but a fully general theory will likely not be straightforward to develop
because of the aforementioned boundary layers.

Based on the oscillatory nature of the cell resonance error, we also 
describe a technique to asymptotically 
reduce it by taking a weighted average at several 
domain sizes. The method is motivated by the form of the one dimensional corrector
\eqref{eq:one_d_corrector}, but it is designed to work more generally in higher
dimensions. It makes use of
averaging kernels that possess regularity properties and vanishing 
``negative'' moment conditions that accelerate the error's convergence to zero
at large domain sizes relative to $\epsilon$. The kernels are 
similar to, but distinct from, those utilized in \cite{tsai:2005,runborg:2016,leitenmaier2022upscaling}.
In dimension $d=1$ we prove the boundary error then
decays as 
\begin{equation}\label{eq:main_1d_result}
\bigoh{\epsilon/\eta}^{\min\{p+1, q+3\}}
\end{equation}
 where $q$ and $p$ depend on the regularity and vanishing moment
properties of the kernel, respectively. 
Since the proof only relies on the 
decomposition \eqref{eq:one_d_corrector}, the result would immediately 
generalize if \eqref{eq:one_d_corrector} were shown to hold for $d> 1$. 
In lieu of a fully general theory, however,
we present numerical 
examples that suggest the strategy indeed is effective in higher dimensions.


As mentioned above, the proposed method does not modify the form of the original cell
problem from homogenization theory. Instead, it involves solving a series of 
microscale problems at different cell sizes $\eta$. Hence, there is no need to develop additional 
approximations to other operators, in contrast to existing methods proposed in the literature. 
In this regard, the method is well suited to be combined with existing reduced basis 
(RB) techniques in the numerical homogenization literature 
\cite{abdulle2012reduced,abdulle2015reduced,boyaval2008reduced,nguyen2008multiscale}, 
where the boundary resonance error is an open problem. Although we do not address 
this in the current work, we offer some further discussion in \cref{sec:discussion} below.

The remainder of the paper is structured as follows.
\Cref{sec:averaging_kernels} below introduces the space 
of averaging kernels $\mathbb{K}^{-p,q}$ and some of its properties. 
Since the results are similar to others available in the literature, 
their proofs can be found in  
\cref{sec:appendix}. In \cref{sec:1d_err_proof} we derive the 
expressions \eqref{eq:cell_res_with_corrector} and 
\eqref{eq:one_d_corrector} and then use the averaging kernels to prove
the cell resonance error can be asymptotically reduced to 
\eqref{eq:main_1d_result}. After deriving an expression for the corrector
$R$ in \eqref{eq:cell_res_with_corrector}
for microscale domains of the form \eqref{eq:I_eta_tubular} in 
\cref{sec:tubular_error}, a numerical strategy for general situations 
is described in \cref{subsec:algo_description}. Numerical 
results are presented and discussed in \cref{sec:numerical_examples}, and then
concluding remarks are offered in \cref{sec:conclusions}. 


\section{Averaging kernels}
\label{sec:averaging_kernels}
Let $C(X)$ denote the space of continuous functions on some $X \subseteq \mathbb{R}$, and 
let $C^k(X)$ be the space of functions that are $k$ times continuously differentiable. 
For compact $X$, let $f \in C_0^k(X)$ if $f$ and its $k$ derivatives vanish at the 
endpoints of $X$. 
\begin{definition}\label{dfn:primitive}
For $g\in C(\mathbb{R})$ let $g^{[0]}(x) = g(x)$, and define the primitive 
\begin{equation*}
g^{[k+1]}(x) := \int_0^{x} g^{[k]}(s)\, ds + c_{k+1}, \qquad k = 0,1,\ldots,
\end{equation*}
where the constant is chosen so that $\int_0^1 g^{[k+1]}(s) \,ds = 0$.
\end{definition}
First, recall the following simple result.
\begin{lemma}\label{lemma:primitives}
Let $g \in C(\mathbb{R})$ be 1-periodic, so that $g(x+1) = g(x)$ $\forall x \in \mathbb{R}$, and let $g$ be mean zero, so that $\int_0^1 g(x)\,dx = 0$. Then the primitive $g^{[1]}(x)$ defined above is also continuous, 1-periodic, and bounded. 
\end{lemma}
The result follows from continuity and the mean-zero assumption; see Appendix \ref{sec:appendix}. 

Consider next the ``incorrect" average of a mean-zero oscillatory function. Let $f \in C(\mathbb{R})$ be 1-periodic, 
and fix $\epsilon > 0$ and $\eta > \epsilon$ with $\eta/\epsilon \notin \mathbb{N}$. 
A simple calculation then shows that
the error in taking the incorrect average of length $\eta$ 
can be bounded as
\begin{equation}
\left| \frac{1}{\eta} \int_{0}^{\eta} f\left(\frac{x}{\epsilon}\right) dx \right| 
= 
 \left| \frac{\epsilon}{\eta} \int_0^{\eta/\epsilon - \floor{\eta/\epsilon}} f(u) du \right|  
\le C\, \norm{f}_{\infty}  \left( \frac{\epsilon}{\eta}\right), \label{eq:arith_averaging_bound}
\end{equation}
where $C = (\eta/\epsilon - \floor{\eta/\epsilon}) \in (0,1)$. Thus, the error decreases like $\epsilon/\eta$. 
This error rate can be asymptotically improved, however, if the average is taken against 
smooth kernels with nice regularity properties.

We next define such averaging kernels that additionally possess vanishing ``negative'' moments and will be
used in the numerical homogenization procedure. These are similar to other kernels with 
vanishing ``positive'' moments found the literature; see e.g. \cite{tsai:2005,runborg:2014,runborg:2016}.
\begin{definition}
Let $p$ and $q$ be nonnegative integers. Define
$K \in \mathbb{K}^{-p, q}\left([1,2]\right)$ if $K \in C^q_0([1,2])$,
$K^{(q+1)} \in BV([1,2])$, 
and 
\begin{equation}\label{eq:vanishing_moment_cond}
\int_1^2 K(x) x^{-r} dx = \begin{cases} 1, \qquad r = 0 \\ 0, \qquad r = 1,\ldots, p. \end{cases}
\end{equation}
\end{definition}
Note that for $K_{\eta}(x):=\frac{1}{\eta}K(x/\eta)$, $K\in\mathbb{K}^{-p,q}\left([1,2]\right) 
\iff K_{\eta} \in \mathbb{K}^{-p,q}\left([\eta,2\eta]\right)$.

\begin{figure}[h]       
    \centering
 \includegraphics[width=0.45\textwidth]{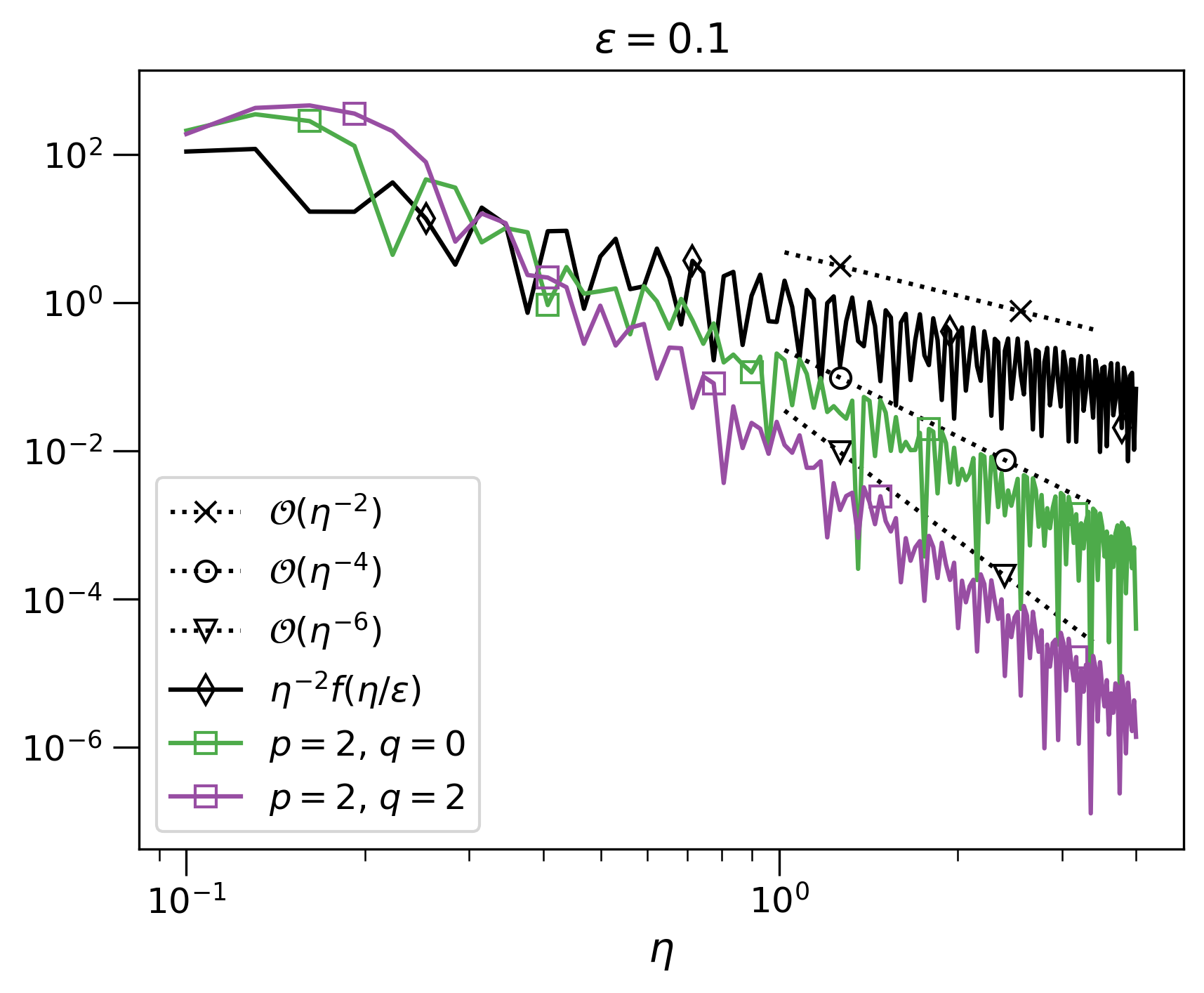}
 \includegraphics[width=0.45\textwidth]{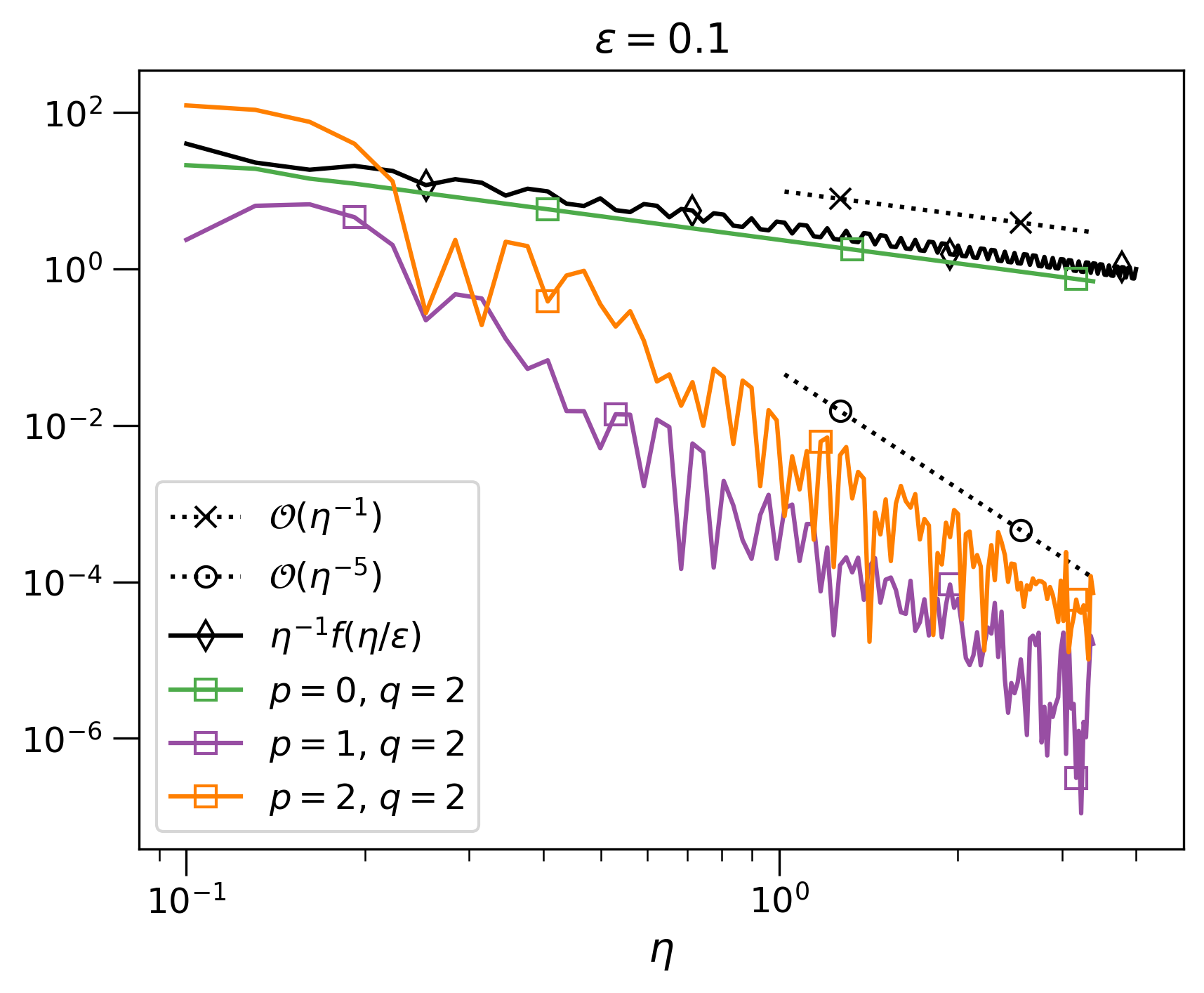}
\caption{ 
The decay of $\eta^{r} f(\eta/\epsilon)$ and $\Upsilon_r(\eta)$ versus $\eta$ for various
averaging kernels $K \in \mathbb{K}^{-p,q}$ of the form \eqref{eq:kernel_form}. 
(Left) $r=2$ and $f(x)=1.1+\sin(2\pi x)$.
(Right) $r=1$ and $f(x)$ is the even periodic extension of the piecewise constant function
$\chi_{[0,1/4]} + \chi_{[1/4,1/2]}$. 
} 
\label{fig:lemma_demonstration}
\end{figure}

The following lemma establishes the bounds needed to prove the main result in 
\cref{sec:1d_err_proof}. 
Since analogous results can be 
found in the literature for different kernels,\footnote{See, for example, 
Lemma 2.3 in \cite{runborg:2014}.} the proof is moved to \cref{sec:appendix}.
\begin{lemma}\label{lem:avg_kernels2}
Let $K \in \mathbb{K}^{-p, q}\left([1,2]\right)$, and let $K_{\eta}(x) := \frac{1}{\eta} K(x/\eta)$ for $\eta >0$. 
Let $f \in C(\mathbb{R})$ be 1-periodic, and let $\chevron{f} = \int_0^1 f(s)ds$.  
Let $0 < \epsilon < \eta $, and for $r \in \mathbb{N}$ let  
\begin{equation*}\label{eq:averaged_decay}
\Upsilon_r(\eta):= \int_{\eta}^{2\eta} K_{\eta}(x) x^{-r} f\left(\frac{x}{\epsilon}\right) \, dx .
\end{equation*}
Then: 
\begin{equation} \label{eq:kernel_decay_estimate}
\left| \Upsilon_r(\eta) \right| \le 
\begin{cases}
C_1\, \norm{f}_{\infty} r^{q+1}\, \left(\frac{\epsilon}{\eta}\right)^{q+2} \eta^{-r}, \qquad r \le p \\
C_1\, \norm{f}_{\infty} r^{q+1} \, \left(\frac{\epsilon}{\eta}\right)^{q+2} \eta^{-r} 
+ \left|\chevron{f}\right| \norm{K}_{\infty} \eta^{-r}, \qquad r > p
\end{cases}
\end{equation}
where $C_1$ is independent of $\epsilon$, $\eta$, and $r$, and depends only on $q$, $K$ and its derivatives.
\end{lemma}
For fixed $\epsilon$, \cref{fig:lemma_demonstration} compares the decay as 
a function of $\eta$ of $\eta^{-r} f(\eta/\epsilon)$ and $\Upsilon_r(\eta)$ for two different periodic 
functions $f$; the rates in \cref{lem:avg_kernels2} are seen to be 
sharp. The averaging kernels $K \in \mathbb{K}^{-p,q}$ used are of the form 
\begin{equation} \label{eq:kernel_form}
 K(x) = (x-1)^{q+1}(x-2)^{q+1} \sum_{j=0}^p a_j x^j,
\end{equation}
where the coefficients $\{a_j\}_{j=0}^p$ are chosen to ensure \eqref{eq:vanishing_moment_cond}
holds. 

Having shown the requisite facts for the function class $\mathbb{K}^{-p,q}$, we are now ready to 
develop our strategy for reducing the boundary error.

\section{Reducing the cell resonance error}   

\subsection{One dimension}
\label{sec:1d_err_proof}
Here and throughout the section, let $a:\mathbb{R}\to\mathbb{R}$ be a 1-periodic, continuous function
that satisfies 
\begin{equation*}
0 < \lambda < a(s) < \Lambda, \qquad \forall s \in \mathbb{R}.
\end{equation*}
In one dimension the cell problem \eqref{eq:cell_problem} rescaled by $\epsilon$ is
\begin{align}\label{eq:1d_true_cell_problem}
-\frac{d}{dx} \Big( a^{\epsilon}(x) &\frac{d\chi}{dx} (x) \Big) = \frac{d}{dx} a^{\epsilon}(x), \qquad x\in [-\epsilon/2,\epsilon/2] \\
\chi \,\,\,\,& \epsilon \text{-periodic and mean zero} \nonumber
\end{align}
with corresponding homogenized coefficient
\begin{equation}\label{eq:true_1d_har_full}
\overline{a}  
=\left(\frac{1}{\epsilon}\int_{-\epsilon/2}^{\epsilon/2} \frac{1}{a^{\epsilon}(s)}ds\right)^{-1}, 
\end{equation}
where $a^{\epsilon}(x) = a(x/\epsilon)$. 
To study the cell resonance error we consider the boundary value problem \eqref{eq:1d_true_cell_problem} posed
instead on the domain $I_{\eta} = [-\eta/2,\eta/2]$, where $\eta$ is inconsistent with the periodicity 
$\epsilon$ of the coefficient $a^{\epsilon}$.
The corresponding ``incorrect'' harmonic average is then 
\begin{equation}\label{eq:incorrect_1d_har_full}
\left(\frac{1}{\eta}\int_{-\eta/2}^{\eta/2} \frac{1}{a^{\epsilon}(s)}ds\right)^{-1}. 
\end{equation}
After rescaling the integral by $\epsilon$, the periodic integrand $1/a(s)$ can be 
decomposed into the sum of odd and even contributions, 
and since the odd part will vanish by symmetry, \eqref{eq:incorrect_1d_har_full} becomes
\begin{equation*}
\left(\frac{2\epsilon}{\eta}\int_{0}^{\eta/2\epsilon} b_{\rm even}(s) ds\right)^{-1}, 
\end{equation*}
where $b_{\rm even}$ denotes the even part of $b(s) := 1/a(s)$. In the development below, 
we hence assume without loss of generality that $b$ is an even function, and for $\eta > \epsilon/2$
we define the ``incorrect'' harmonic average  
\begin{equation}\label{eq:naive_avg}
\rho^{\epsilon} (\eta) := 
\left(\frac{\epsilon}{\eta}\int_0^{\eta/\epsilon} b(s) \, ds\right)^{-1}.
\end{equation}
%
After a simple algebraic rearrangement, a calculation similar to \eqref{eq:arith_averaging_bound} gives
\begin{equation*}\label{eq:one1_cell_res_error}
\left| \overline{a} - \rho^{\epsilon}(\eta)\right| \le C \left(\frac{\epsilon}{\eta}\right),
\end{equation*}
which is the one-dimensional version of the cell-resonance error first derived in \cite{pingbingming:2004}.
Within this bound, however, the error has a richer structure, which we show next. 
\begin{lemma}\label{lem:expansion}
Let
 \begin{equation}\label{eq:B_defn}
B = \norm{\frac{1}{a} - \int_0^1 \frac{1}{a(s)} ds}_{L^{\infty}(\mathbb{R})}
\end{equation}
be the supremum of the fluctuations about the mean for the function $1/a$, and let
$\overline{a}$ be the true harmonic average \eqref{eq:true_1d_har_full}. If $\eta/\epsilon > \overline{a}B$, 
then $\rhoeps$ is given by the expansion
\begin{equation}\label{eq:main_decomposition}
\rhoeps(\eta) = \overline{a} \sum_{k=0}^{\infty} \left(-\frac{\epsilon}{\eta} \, q(\eta/\epsilon) \right)^k
\end{equation}
where $q(\eta/\epsilon)$ is an $\epsilon$-periodic function of $\eta$. 
\end{lemma}
\textit{Proof:} 
Let $b(s):= 1/a(s)$, and let $\chevron{b} := \int_0^1 b(s) \, ds$ be its 
arithmetic average. Note $\overline{a} = 1/\chevron{b}$. Then
\begin{align*}
\rho^{\epsilon}(\eta) 
&= \left(\chevron{b} +  \frac{\epsilon}{\eta} \int_0^{\eta/\epsilon} \left(b(s)-\chevron{b}\right) \, ds \right)^{-1} \\ 
&= \overline{a} \left(1  + \overline{a} \,\,\frac{\epsilon}{\eta} \int_0^{\eta/\epsilon} \left(b(s)-\chevron{b}\right) \, ds \right)^{-1} \\
&= \overline{a} \left(1 + z\left(\eta/\epsilon\right)\right)^{-1},
\end{align*}
where
$$z(\eta/\epsilon) = \overline{a}\,\, \frac{\epsilon}{\eta} \int_0^{\eta/\epsilon} \left(b(s)-\chevron{b}\right) ds.$$
The condition $\eta/\epsilon > \overline{a} B$ ensures that $|z(\eta/\epsilon)| < 1$ and hence that $\rhoeps$ can 
be expanded in a uniformly convergent geometric series, as desired. Since $b(s)-\chevron{b}$ is a mean-zero, periodic
function, 
\begin{equation*}\label{eq:q_defn}
 q(\eta/\epsilon) := z(\eta/\epsilon) \frac{\eta}{\epsilon} = \overline{a} \int_{0}^{\eta/\epsilon}( b(s)-\chevron{b} ) ds 
\end{equation*}
is also $\epsilon$-periodic by \cref{lemma:primitives}. $\square$

The decomposition of the cell resonance error \eqref{eq:main_decomposition} motivates the use
of weighted averages against kernels $K \in \mathbb{K}^{-p,q}$, which we define next. 
\begin{definition}\label{dfn:smooth_avg}
For $K\in \mathbb{K}^{-p,q}\left([1,2]\right)$, $\eta > 0$ and $K_{\eta}(x) = \frac{1}{\eta}K(x/\eta)$, 
define 
\begin{equation}\label{eq:smooth_avg}
S^{\epsilon}(\eta) := \int_{\eta}^{2\eta} K_{\eta}(x)\rho^{\epsilon}(x)\, dx.
\end{equation}
\end{definition}
The averaged error \eqref{eq:smooth_avg} accelerates the convergence of $\rhoeps(\eta)$ 
to the true harmonic average $\overline{a}$ at a rate that depends on $p$ and $q$, which we
now show. 
\begin{theorem}\label{thm:1d_err_thm}
For non-negative integers $p$ and $q$ let $K \in \mathbb{K}^{-p,q}\left([1,2]\right)$, and let $\epsilon, \eta > 0$ such that 
\begin{equation}\label{eq:converge_cond1d}
\overline{a} B\left(1+\frac{1}{p+1}\right)^{q+1} <  \frac{\eta}{\epsilon}
\end{equation}
where $B$ is given by \eqref{eq:B_defn} and
$\overline{a}$ is the true harmonic average \eqref{eq:true_1d_har_full}. 
Let $S^{\epsilon}$ be defined by \eqref{eq:smooth_avg}.
Then 
\begin{equation}\label{eq:1d_err_estimate}
\left|\, S^{\epsilon}(\eta) - \overline{a} \, \right| \le   
C \, \varphi\big(\overline{a}B \frac{\epsilon}{\eta}\big)\left(\frac{\epsilon}{\eta}\right)^{\min\{p+1,q+3\}},
\end{equation}
where $C$ depends on $p$, $q$, $K$, $\overline{a}$ and $B$, and $\varphi(\theta) := 1/(1-\theta)$. 
\end{theorem}
\textit{Proof:}
Since the series expansion \eqref{eq:main_decomposition} converges uniformly and absolutely, it can 
be inserted into the weighted average \eqref{eq:smooth_avg}, and order of the sum 
and integral can be interchanged. The resulting expression is
\begin{align*}
S^{\epsilon}(\eta) &= \int_{\eta}^{2\eta} K_{\eta}(x) \rho^{\epsilon}(x) dx  \\
&= \overline{a} + \sum_{r = 1}^{\infty} (-1)^r \overline{a}^{r+1}  \epsilon^r 
\int_{\eta}^{2 \eta} \Big[K_{\eta}(x) x^{-r} 
\overbrace{\left( \int_0^{x/\epsilon} \left(b(s)-\chevron{b}\right)ds\right)^r \, }^{G_r(x/\epsilon) := }\, \Big] dx, 
\end{align*}
which implies 
\begin{equation}\label{eq:1d_err_inf_sum}
\left| S^{\epsilon}(\eta) - \overline{a}\right| \le 
\sum_{r=1}^{\infty} \overbrace{\overline{a}^{r+1} \epsilon^r 
\left|\int_{\eta}^{2 \eta} K_{\eta}(x) x^{-r} G_r(x/\epsilon) dx \right|}^{F_r:=}.
\end{equation}
Fixing $r$, each individual term in the summations can be bounded using \cref{lem:avg_kernels2}. Since 
$b-\chevron{b}$ periodic, and mean-zero, and continuous, \cref{lemma:primitives} implies $G_r$ is both periodic and continuous as well. 
Hence, the assumptions of \cref{lem:avg_kernels2} are indeed satisfied, so that  
\begin{align}
F_r 
&\le \overline{a}^{r+1} \epsilon^r \cdot \begin{cases} C(q,K) \norm{G_r}_{\infty} \,r^{q+1} \left(\frac{\epsilon}{\eta}\right)^{q+2} \eta^{-r}, \qquad r \le p \\ 
C(q,K) \norm{G_r}_{\infty} \, r^{q+1}\left(\frac{\epsilon}{\eta}\right)^{q+2} \eta^{-r} + \left|\chevron{G_r}\right|\norm{K}_{\infty} \eta^{-r}, \qquad r > p \end{cases} \nonumber \\
&= \begin{cases} C(q,K) \norm{G_r}_{\infty} r^{q+1}\, \overline{a}^{r+1} \left(\frac{\epsilon}{\eta}\right)^{r+q+2}, \qquad r \le p  \\ 
C(q,K) \norm{G_r}_{\infty}\, r^{q+1}\, \overline{a}^{r+1} \left(\frac{\epsilon}{\eta}\right)^{r+q+2} 
+ \left|\chevron{G_r}\right| \norm{K}_{\infty}\, \overline{a}^{r+1} \left(\frac{\epsilon}{\eta}\right)^{r}, \qquad r > p.\end{cases}\label{eq:Fr_bound}
\end{align}
Furthermore, because $b-\chevron{b}$ is periodic and mean zero, 
\begin{equation}\label{eq:bound_antideriv_of_b}
\sup_{x\in\mathbb{R}} \left| \int_0^{x/\epsilon} (b(s)-\chevron{b})\, ds \right| \le B,
\end{equation}
and hence
\begin{equation}\label{eq:Gr_bound}
\norm{G_r}_{\infty} \le B^r; 
\end{equation}
an identical bound exists for $|\chevron{G_r}|$. 
The infinite sum \eqref{eq:1d_err_inf_sum} can be broken up into two pieces, at which
point the bounds \eqref{eq:Fr_bound} and \eqref{eq:Gr_bound} can be inserted to 
bound the difference between the true harmonic average and $S^{\epsilon}(\eta)$ by
three parts:
\begin{align}
\left| S^{\epsilon}(\eta) - \overline{a}\right| \le  
\sum_{r=1}^p &C(q,K) r^{q+1} \overline{a}^{r+1} B^r \left(\frac{\epsilon}{\eta}\right)^{r+q+2} \nonumber \\
+ & \sum_{r=p+1}^{\infty}\Big[ C(q,K) r^{q+1} \overline{a}^{r+1} B^r \left(\frac{\epsilon}{\eta}\right)^{r+q+2} 
+  \norm{K}_{\infty}\overline{a}^{r+1} B^r \left(\frac{\epsilon}{\eta}\right)^{r} \Big],  \label{eq:rhs_of_inequality}
\end{align}
The first part on the right hand side (RHS) 
of the inequality \eqref{eq:rhs_of_inequality} is a finite sum, which can be bounded as
\begin{equation}\label{eq:first_bound}
\sum_{r=1}^p C(q,K) r^{q+1} \overline{a}^{r+1} B^r \left(\frac{\epsilon}{\eta}\right)^{r+q+2}
\le  C(q,K) \, \overline{a}^2B\, p^{q+2}\left(\frac{\epsilon}{\eta}\right)^{q+3}.
\end{equation}
The third part of the RHS of \eqref{eq:rhs_of_inequality} is the tail of a geometric series. The assumption 
\eqref{eq:converge_cond1d} implies 
\begin{equation}
\sum_{r=p+1}^{\infty}\norm{K}_{\infty}\overline{a}\, (\overline{a} B)^r \left(\frac{\epsilon}{\eta}\right)^{r}  
= \norm{K}_{\infty}\overline{a}\,(\overline{a} B)^{p+1}\left(\frac{\epsilon}{\eta}\right)^{p+1} \varphi\big(\overline{a}B \frac{\epsilon}{\eta}\big).
\label{eq:second_bound}
\end{equation}
Finally, the second part of the RHS of \eqref{eq:rhs_of_inequality} can be bounded with the ``integral test" 
from elementary calculus, which says
\begin{equation*}
\sum_{r=p+1}^{\infty}f(r) \le f(p+1) + \int_{p+1}^{\infty} f(r)\,dr,
\end{equation*}
provided $f$ is a monotonically decreasing function of $r$. Here $f(r)=r^{q+1}(\overline{a} B)^r\left(\frac{\epsilon}{\eta}\right)^{r}$, 
and the assumption \eqref{eq:converge_cond1d} guarantees that $f(r+1) < f(r)$ indeed holds for every $r>p$. 
Hence
\begin{align}
\sum_{r=p+1}^{\infty} C(q,K) r^{q+1} (\overline{a} B)^r \left(\frac{\epsilon}{\eta}\right)^{r+q+2}
 \le C(q,K)  &\Big[ 
(p+1)^{q+1}(\overline{a} B)^{p+1}\left(\frac{\epsilon}{\eta}\right)^{p+q+3} \nonumber \\
&+ \left(\frac{\epsilon}{\eta}\right)^{q+2}\underbrace{\int_{p+1}^{\infty} r^{q+1}\alpha^r \, dr}_{(\ast)}\Big],   \label{eq:second_piece_calc}
\end{align}
where $\alpha = \overline{a}B\epsilon/\eta$. The integral $(\ast)$ can be computed with $q+1$ applications
of integration by parts:
\begin{align}
(\ast) &= -\int_{p+1}^{\infty} (q+1)r^q\left(\log(\alpha)\right)^{-1} \alpha^r dr + (p+1)^{q+1}\left(-\log(\alpha)\right)^{-1} \alpha^{p+1} \nonumber \\
&= \int_{p+1}^{\infty} (q+1)!\left(-\log(\alpha)\right)^{-(q+1)}\alpha^r dr \nonumber \\
&+ \sum_{k=1}^{q+1}\frac{(q+1)!}{(q+2-k)!}(p+1)^{q+2-k} \left(-\log(\alpha)\right)^{-k}\alpha^{p+1}\nonumber \\
&= (q+1)!\left(-\log(\alpha)\right)^{-(q+2)}\alpha^{p+1}  \label{eq:third_bound} \\
&+ \alpha^{p+1}\sum_{k=1}^{q+1}\frac{(q+1)!}{(q+2-k)!}(p+1)^{q+2-k} \left(-\log(\alpha)\right)^{-k}\nonumber
\end{align}
Comparing the dependence on $\epsilon/\eta$ in the upper bounds on the three parts on the RHS of the inequality 
\eqref{eq:rhs_of_inequality} (given in \eqref{eq:first_bound}, \eqref{eq:second_bound}, and both \eqref{eq:second_piece_calc} and \eqref{eq:third_bound}) then gives the desired result:        
\begin{equation*}
\left| S^{\epsilon}(\eta) - \overline{a} \right| \le C(p,q,K,\overline{a}, B) \varphi\big(\overline{a}B \frac{\epsilon}{\eta}\big)
\left(\frac{\epsilon}{\eta}\right)^{\min\{p+1,q+3\}} \square
\end{equation*} 

Next we note that in the special case that $K(x) = 1 \in \mathbb{K}^{0,-1}$, \eqref{eq:smooth_avg} becomes a simple arithmetic average
of $\rhoeps$, and the result \eqref{eq:1d_err_estimate} implies that the averaged error asymptotically decays as $\epsilon/\eta$. 
A more direct calculation, however, reveals that this estimate can be improved to second order in $\epsilon/\eta$. 
\begin{theorem}\label{thm:second_order}
Let $\eta/\epsilon > \overline{a}B$. Then 
\begin{equation}\label{eq:simple_avg_result}
\left|\frac{1}{\eta} \int_{\eta}^{2\eta} (\rhoeps(x) - \overline{a})\, dx \right| \le
C \, \varphi\big(\overline{a}B \frac{\epsilon}{\eta}\big)\left(\frac{\epsilon}{\eta}\right)^2
\end{equation}
where $C$ depends on $\overline{a}$ and $B$, and $\varphi(\theta) := 1/(1-\theta)$. 
\end{theorem}
Since the proof stategy is similar to that of \cref{thm:1d_err_thm}, it is presented in 
\cref{sec:appendix_second_order_decay}. 

In practice, we find in the numerical experiments presented in \cref{sec:numerical_examples} that 
$\norm{K}_{\infty}$ grows relatively fast as $p$ and $q$ grow, and hence
the error decay rate predicted by \eqref{eq:1d_err_estimate} is not seen until $\eta/\epsilon$
is relatively large. Although the decay rate \eqref{eq:simple_avg_result} is only second order in $\epsilon/\eta$,
this relatively simple 
averaging strategy is found to be most effective in practice, particularly in higher dimensions 
where solving the cell problem \eqref{eq:cell_problem} on a domain with size $\delta := \eta/\epsilon$ 
in each direction 
becomes computationally expensive at large $\delta$. 


\subsection{Tubular domains in two dimensions}
\label{sec:tubular_error}
Next, 
let dimension $d=2$, and in the following assume that $a(x_1,x_2)$ is an isotropic, coercive tensor that 
is $1$-periodic in both arguments. Because we make use of classical Sturm-Liouville theory below,
we additionally assume that $a$ is once continuously differentiable. 

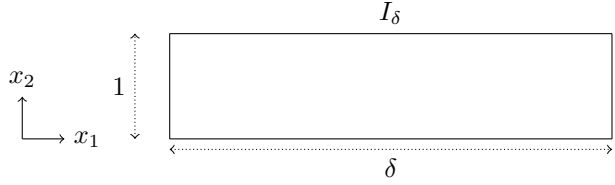
\begin{figure}
\centering
\begin{tikzpicture}[xscale=1.4,yscale=1.4]
\draw[densely dotted] (2.10,1)node[above]{$I_{\delta}$};
\draw[densely dotted, <->] (-0.33,0) -- (-0.33,1);
\draw (-0.33,0.5)node[left]{$1$}; 
\draw[->](-1.4,0)--(-1.4,0.4)node[above]{$x_2$}; 
\draw[->](-1.4,0)--(-1.0,0)node[right]{$x_1$}; 
\draw[](0,0) -- (4.2,0);
\draw[](0,1) -- (4.2,1);
\draw[](0,1) -- (0,0);
\draw[](4.2,1) -- (4.20,0);
\draw[densely dotted, <->] (0,-0.1) -- (4.2,-0.1); 
\draw (2.1,-0.1)node[below]{$\delta$}; 
\end{tikzpicture}
\caption{A two-dimensional ``tubular" domain $I_{\delta} = [-\delta/2,\delta/2]\times[-1/2,1/2]$,
$\delta \notin \mathbb{N}$.
}
\label{fig:domain_cartoon}
\end{figure}

For $\delta = \eta/\epsilon$, consider domains of the form 
\begin{equation*}
I_{\delta} = [-\delta/2,\delta/2]\times[-1/2,-1/2],
\end{equation*}
as depicted in \cref{fig:domain_cartoon}. Assume $\delta \notin\mathbb{N}$, so that 
there is a mismatch between the size of the micro domain and the periodicity of $a$
in the $x_1$ direction.  
We show that the cell resonance error can be written in the form \eqref{eq:cell_res_with_corrector} 
from the introduction, where the ``corrector'' $R(\delta)$ is given by the average
of a locally $\epsilon$-periodic function.

Consider first the periodic extension to all of $\mathbb{R}^2$
of the solution $\chi_j$, $j\in\{1,2\}$, to the true 
cell problem \eqref{eq:cell_problem} (i.e.\ with $\delta =1$).
Next consider the solution $\phi_{\delta,j}$ to the cell problem posed on 
the ``incorrect" domain $I_{\delta}$: 
\begin{align}
-\nabla \cdot \big(a \nabla \phi_{\delta,j}\big) &= \nabla \cdot \big( ae_j\big), \qquad I_{\delta} \nonumber \\
\int_{I_{\delta}} \phi_{\delta,j} \, dx &= 0   \label{eq:phi_eta_j_problem}
\end{align}
with periodic boundary conditions, and define the difference 
\begin{equation}\label{eq:u_defn}
u_{\delta,j}(x_1, x_2) := \phi_{\delta,j}(x_1,x_2) - \chi_j(x_1, x_2). 
\end{equation}
Note that if $\delta \in \mathbb{N}$, then $u_{\delta,j} =0$ identically. 
The standard elliptic theory guarantees that both $\chi_j$ and $\phi_{\delta,j}$,
(and hence $u_{\delta,j}$) are twice continuously differentiable. 
More generally $u_{\delta,j}$ solves the following elliptic problem in $I_{\delta}$: 
\begin{align}
-\nabla \cdot \big(a \nabla u_{\delta,j}\big) = 0&, \qquad I_{\delta} \nonumber  \\
u_{\delta,j}(x_1,-1/2) = u_{\delta,j}(x_1,1/2)&, \qquad x_1 \in [-\delta/2,\delta/2] \nonumber \\
u_{\delta,j} = g_{\delta,j}(x_2)&, \qquad x_1=-\delta/2 \nonumber \\ 
u_{\delta,j} = h_{\delta,j}(x_2)&, \qquad x_1 = \delta/2 \label{eq:u_j_eqn_full}
\end{align}
where $g_{\delta, j}$ and $h_{\delta,j}$ are simply the difference between
 $\phi_{\delta,j}$ and $\chi_j$ at the points $x_1=\mp \delta/2$. Since $\delta$ is 
assumed to be inconsistent with the periodicity of $a$, both $g_j$ and $h_j$
 are oscillatory functions of $\delta$ for fixed $x_2$. 
 Both $u_{\delta,j}$ and $\phi_{\delta,j}$
are of course $1$-periodic in $x_2$; we next show that they are also both locally 
$1$-periodic in $x_1$.

\begin{theorem}\label{thm:u_loc_per}
The solution $\phi_{\delta,j}(x_1,x_2)$ to the elliptic problem \eqref{eq:phi_eta_j_problem} 
is locally $1$-periodic in $x_1$. 
\end{theorem}
\textit{Proof:}
Because $\chi_j$ is $1$-periodic, it is sufficient to show that the difference 
$u_{\delta,j} = \phi_{\delta,j} -\chi$ is locally $1$-periodic. Since the proof proceeds
identically for $j \in \{1,2\}$, we drop the subscript below for notational convenience. 

By assumption of isotropy, the governing equation for $u_{\delta}$ simplifies to 
\begin{equation}\label{eq:u_j_eqn}
-\left( \pderiv{}{x_1}\Big( a(x_1,x_2) \pderiv{}{x_1}\Big) + \pderiv{}{x_2} \Big( a(x_1,x_2) \pderiv{}{x_2}\Big)\right)
u_{\delta}(x_1,x_2) = 0.
\end{equation} 
For fixed $x_1 \in \mathbb{R}$, define the symmetric, positive-definite  
differential operator $L(x_1)$ by:
\begin{equation*}
L(x_1)f:= -\pderiv{}{x_2} \left(a(x_1,x_2) \pderiv{}{x_2}f\right), 
\end{equation*}
acting on $C^2$ functions $f: [-1/2,1/2] \to \mathbb{R}$ with 
periodic boundary conditions. 
Sturm-Liouville theory \cite{math_phys_book} 
guarantees that $L(x_1)$ has a countably 
infinite set of eigenvalues $\{\lambda_i(x_1)\}_{i=0}^{\infty}$
that satisfy $\lambda_0(x_1) = 0$ and $\lambda_{i+1}(x_1) > \lambda_i(x_1)$ for 
all $i \ge 0$, as well as an associated set of complete, mutually orthonormal
eigenfunctions $\{\varphi_i(x_1,x_2)\}_{i=0}^{\infty}$. Note also that by 
periodicity of $a$, $L(x_1+1) = L(x_1)$ for any $x_1 \in\mathbb{R}$. Since the
eigenvalues and their orthonormal eigenfunctions are unique, 
this implies that both are $1$-periodic in $x_1$. 

Consider an eigenfunction expansion of the solution to the boundary value problem (BVP) \eqref{eq:u_j_eqn_full} 
\begin{equation}\label{eq:eigen_expansion}
u_{\delta}(x_1,x_2) = \sum_{i=0}^{\infty} \alpha_i(x_1) \varphi_i(x_1,x_2)
\end{equation}
where 
\begin{equation*}
\alpha_i(x_1) = \int_{-1/2}^{1/2} u_{\delta}(x_1,x_2) \varphi_i(x_1,x_2) \, dx_2.
\end{equation*}
Since $u_{\delta}$ is $C^2$, the expansion converges uniformly and absolutely.

Inserting \eqref{eq:eigen_expansion} into \eqref{eq:u_j_eqn} then gives
\begin{equation*}
\sum_{i=0}^{\infty}\Big[ -\pderiv{}{x_1}\Big( a(x_1,x_2) \pderiv{}{x_1}\big( \alpha_i(x_1)\varphi_i(x_1,x_2)\big)\Big)
+\lambda_i(x_1) \alpha_i(x_1) \varphi_i(x_1,x_2) \Big] = 0 
\end{equation*}
Let $\psi_i(x_1,x_2) := \alpha_i(x_1)\varphi_i(x_1,x_2)$. A sufficient
condition for $u$ to solve \eqref{eq:u_j_eqn} is for each $\psi_i$ 
to satisfy the following boundary value problem (BVP) 
 in $x_1$ that is parameterized by $x_2 \in [-1/2,1/2]$:
\begin{align}
-\pderiv{}{x_1}\Big( a(x_1,x_2) \pderiv{}{x_1}\big(\psi_i(x_1,x_2)\big)\Big) \nonumber
+\lambda_i(x_1)\psi_i(x_1,x_2) = 0&, \qquad x_1 \in [-\delta/2,\delta/2] \nonumber\\
\psi_i = \beta_i(-\delta/2) \varphi_i(-\delta/2,x_2)&, \qquad x_1 = -\delta/2 \nonumber  \\
\psi_i = \gamma_i(\delta/2) \varphi_i(\delta/2,x_2)&, \qquad x_1 = \delta/2. \label{eq:psi_bvp}
\end{align} 
Here $\beta_i(-\delta/2)$ is the $i$-th coefficient in the expansion in eigenfunctions of 
$L(-\delta/2)$ of the boundary condition $g_{\delta}(x_2)$ in \eqref{eq:u_j_eqn_full}, 
while $\gamma_i(\delta/2)$ is the coefficient of the expansion of $h_{\delta}(x_2)$ in eigenfunctions of $L(\delta/2)$. 
Since $\lambda_i(x_1) \ge 0$ for all $i \ge 0$ 
and $\forall x_1 \in \mathbb{R}$, \eqref{eq:psi_bvp}
has a unique solution $\forall x_2 \in [-1/2,1/2]$ and $i\ge0$. 

For each solution $\psi_i(x_1,x_2)$, define the ``flux'' 
$$v_i(x_1,x_2):= a(x_1,x_2)\pderiv{}{x_1} \psi_i(x_1,x_2) $$
and consider a reformulation of the BVP \eqref{eq:psi_bvp} into a linear initial value
problem system 
\begin{equation*}\label{eq:ivp_system}
\pderiv{}{x_1} 
\begin{pmatrix}
 v_i(x_1,x_2) \\ 
\psi_i (x_1,x_2)
\end{pmatrix}
= 
\begin{pmatrix}
 0 & \lambda_i(x_1,x_2) \\ 
 1/a(x_1,x_2) & 0 \\
\end{pmatrix}
\begin{pmatrix}
 v_i(x_1,x_2) \\ 
\psi_i(x_1,x_2) 
\end{pmatrix} 
\end{equation*}
for $x_1 > -\delta/2$,
where the initial conditions are simply the values of $\psi_i$ and $v_i$ at $x_1 = -\delta/2$. 
Since both $1/a$ and $\lambda_i$ are $1$-periodic functions of $x_1$, 
Floquet's theorem \cite{floquet1883equations} implies that the solution can be written as
\begin{equation*}
\begin{pmatrix}
 v_i (x_1,x_2)\\ 
\psi_i (x_1,x_2)
\end{pmatrix}
= P_i(x_1,x_2) e^{x_1 K_i(x_2)}  c_i
\end{equation*}
for some $c_i \in \mathbb{R}^2$. Here $K_i$ and $P_i$ are both $2\times 2$
matrices that generally depend on $x_2$; however, $K_i$ is constant in $x_1$ and 
$P_i$ is $1$-periodic in $x_1$. Hence, each $\psi_i$ is 
locally $1$-periodic in $x_1$. Since the 
eigenfunction expansion \eqref{eq:eigen_expansion} 
converges uniformly and absolutely, the infinite sum 
$u_{\delta} = \sum_{i=0}^{\infty} \psi_i$ 
is also locally $1$-periodic in $x_1$, as desired. $\square$
\begin{remark}
Note that while \cref{thm:u_loc_per} characterizes the behavior of 
$u_{\delta, j}(x_1,x_2)$ (and hence $\phi_{\delta,j}(x_1,x_2)$) as a function of 
$x_1$ for \emph{fixed} $\delta$, it does not describe the functional dependence on $\delta$. In general 
the dependence is felt through the oscillatory boundary conditions $g_{\delta, j}$
and $h_{\delta, j}$ in \eqref{eq:u_j_eqn_full}. 
\end{remark}

Next, define 
$$
\rho(\delta):= \intdelta (a + a\nabla \phi_{\delta, j})  \, dx_2 dx_1 
$$
to be the approximation to the homogenized tensor with cell size $\delta$; 
the true homogenized coefficients $\overline{a}$ result when $\delta = 1$, 
or more generally when $\delta \in \mathbb{N}$. As a corollary to  
\cref{thm:u_loc_per}, we next show the cell 
resonance error can be written in the form \eqref{eq:cell_res_with_corrector} 
from the introduction.
\begin{corollary}
Let $E(\delta) = \rho(\delta) - \overline{a}$ be the cell resonance error.
Then 
\begin{equation*}\label{eq:defn_res_error_explicit}
E(\delta) = \frac{1}{\delta} P(\delta) + 
\frac{1}{\delta} \int_{-\delta/2}^{\delta/2} \zeta(x_1; \delta) \, dx_1
\end{equation*}
where $P$ is an isotropic tensor that is $1$-periodic in $\delta$ and 
$\zeta$ is a diagonal tensor that is locally $1$-periodic in $x_1$. 
\end{corollary}
\textit{Proof:}        
First note that since $a$ is isotropic, both $\overline{a}$ and
$\rho(\delta)$ are diagonal: $\overline{a}_{ij} = \rho_{ij}(\delta) = 0$ 
for $i \ne j$. 
Next, define for $\nu \in \{1,2\}$ the arithmetic averages
$$
\chevron{a_{\nu\nu}} = \int_{-1/2}^{1/2}\int_{-1/2}^{1/2} a_{\nu\nu} \, dx_2 dx_1
$$
and 
$$
\chevron{a_{\nu\nu}\pderiv{\chi_{\nu}}{x_{\nu}}} = \int_{-1/2}^{1/2}\int_{-1/2}^{1/2} a_{\nu\nu}\pderiv{\chi_{\nu}}{x_{\nu}} \, dx_2 dx_1
$$
(here Einstein summation is not implied), 
so that $\overline{a}_{\nu\nu} = \chevron{a_{\nu\nu}} + \chevron{a_{\nu\nu} \partial \chi_{\nu}/\partial x_{\nu}}$.
Using $\phi_{\delta,\nu} = u_{\delta,\nu} + \chi_{\nu}$ from \eqref{eq:u_defn},  
the resonance error is given by 
\begin{equation}\label{eq:decompose_error}   
E_{\nu\nu}(\delta)= \intdelta \Big( 
p_{\nu}(x_1,x_2)
+ a_{\nu\nu}\pderiv{u_{\delta,\nu}}{x_{\nu}}(x_1,x_2) 
\Big) dx_2 dx_1
\end{equation}
where 
$$p_{\nu}(x_1,x_2)) := \big(a_{\nu\nu}-\chevron{a_{\nu\nu}}\big) + 
\Big(a_{\nu\nu}\pderiv{\chi_{\nu}}{x_{\nu}} - \chevron{a_{\nu\nu}\pderiv{\chi_{\nu}}{x_{\nu}}} \Big)$$
is a mean-zero, $1$-periodic function of both $x_1$ and $x_2$. 
Let
$$ \widetilde p_{\nu}(x_1):=  \int_{-1/2}^{1/2} p_{\nu}(x_1,x_2) \, dx_2, $$
and break $\widetilde p_{\nu}$ into even and odd parts. Then the first part of 
\eqref{eq:decompose_error} is
$$
\intdelta p_{\nu}(x_1,x_2) \, dx_2 dx_1 =   \frac{1}{\delta} \, 2\int^{\delta/2}_0 \widetilde p_{\nu, \rm even}(x_1) \, dx_1 
=: \frac{1}{\delta} P(\delta),
$$
and \cref{lemma:primitives} guarantees that the primitive $P(\delta)$ is $1$-periodic.

The second term in the resonance error \eqref{eq:decompose_error} is 
$$
R(\delta):= \intdelta a_{\nu\nu}(x_1,x_2) \pderiv{u_{\delta,\nu}}{x_{\nu}}(x_1,x_2) \, dx_2dx_1.
$$
By \cref{thm:u_loc_per}, $u_{\nu}$ is 
locally $1$-periodic in $x_1$ for each $\delta$. Since the gradient $\nabla u_{\nu}$ is also 
locally $1$-periodic, setting 
$$
\zeta_{\nu\nu}(x_1;\delta) := \int_{-1/2}^{1/2} a_{\nu\nu}(x_1,x_2)\pderiv{u_{\delta,\nu}}{x_{\nu}}(x_1,x_2) \, dx_2
$$
for $\nu \in \{1,2\}$ gives the desired result. 
$\square$              

\subsection{General numerical strategy}
\label{subsec:algo_description}
%
%

Motivated by the idea that the cell resonance error is an oscillatory 
function of domain size, we now propose a numerical strategy for 
reducing the boundary error in general, higher dimensional settings.


Let the dimension $d > 1$, and for $\eta>0$, let $I_{\eta} = [-\eta/2,\eta/2]^d$. Let 
$\aeps: \mathbb{R}^d \to \mathbb{R}^d$ to be a bounded, coercive function with entries 
that are supported at frequencies $\sim 1/\epsilon$. In general $\aeps$ need not be
periodic. For $1 \le j \le d$, let $\chi_{\eta,j}$ be the solution to the elliptic problem 
\begin{align}
-\nabla \cdot \Big( \aeps(x) \nabla \chi_{\eta,j}(x)\Big) &= \nabla \cdot (\aeps(x) e_j), 
\qquad x \in I_{\eta} \nonumber  \\
\int_{I_{\eta}} \chi_{\eta,j}  \, dx &= 0   \label{eq:chi_j_problem}
\end{align}
posed with periodic boundary conditions. Other boundary conditions such as Dirichlet or 
Neumann are possible; however, from \cite{yue:2006} it is known that periodic conditions 
introduce the least boundary effects. 
In general we assume \eqref{eq:chi_j_problem} 
possesses a suitably defined weak solution, but since we are ultimately interested 
in numerical solutions, we do not further consider any issues of regularity. 

Define for $1 \le i,j \le d$ the estimate for 
the homogenized tensor
\begin{equation}  \label{eq:general_rho_ij}
\rho_{ij}(\eta):= \frac{1}{|I_{\eta}|} \int_{I_{\eta}} \Big(a_{ij}(x) + a_{ik} \pderiv{\chi_{\eta,j}}{x_k}(x) \Big) dx ,
\end{equation}
where summation over $k$ is implied.
Note that if $\aeps$ is $\epsilon$-periodic, for example, the true homogenized coefficient $\overline{a}_{ij}$ is
simply given by \eqref{eq:general_rho_ij} for $\eta = \epsilon$, while 
if $\aeps$ is quasi-periodic or stationary ergodic, then $\overline{a}_{ij}$ is 
given by the limit as $\eta\to\infty$ of \eqref{eq:general_rho_ij},  
where $\chi_{\eta,j}$ is replaced by a function $\chi_{\infty,j}$ that solves the elliptic problem \eqref{eq:chi_j_problem}
posed on $\mathbb{R}^d$ \cite{gloria:2011,kozlov:1994}.

To improve the approximation of the true homogenized coefficients, we introduce the averaged coefficients 
\begin{equation} \label{eq:S_ij}
S_{ij}(\eta) := \int_{\eta}^{2\eta} \Keta(t) \rho_{ij}(t)\, dt, 
\end{equation}
where $\Keta \in \mathbb{K}^{-p,q}([\eta,2\eta])$ for some nonnegative integers $p$ and $q$.
In practice, this integral is approximated by some quadrature rule defined by $M$ weights $\{\omega_m\}_{m=1}^M$
and nodes $T_M := \{t_m\}_{m=1}^M \subset [\eta,2\eta]$. 

With these preliminary definitions, we now summarize how the smoothed approximation 
\eqref{eq:S_ij} of the true homogenized coefficients are obtained in practice. 

\smallskip
\begin{enumerate}
\item[(1)]
For each $t_m \in T_M$: 
\smallskip
   \begin{enumerate}
   \item[(i)] Solve the elliptic problem \eqref{eq:chi_j_problem} with $\eta = t_m$ for each $1 \le j \le d$
   using e.g.\ a finite element method, finite difference method, etc.
\smallskip
   \item[(ii)] Using a quadrature method on $I_{t_m}$ and a suitable estimate for the gradient of 
   $\chi_{t_m,j}$, numerically compute and store $\rho_{ij}(t_m)$ for each $1 \le i,j \le d$.
   \end{enumerate}

\smallskip
\item[(2)]
Using the results, estimate the averaged coefficients $S_{ij}$ with
\begin{equation} \label{eq:final_quad}
\sum_{m=1}^M \omega_m \Keta(t_m) \rho_{ij}(t_m). 
\end{equation}
\end{enumerate}
Numerical examples in two dimensions are described below in Section \ref{subsec:results2d}.

\section{Numerical examples}
\label{sec:numerical_examples}

\subsection{Dimension $d=1$}
\label{sec:numerical_results1d}
We now present some one dimensional numerical examples of the proposed 
averaging method for reducing the cell resonance error. 
Since the error can be computed either analytically or with 
quadrature formulas, it is relatively easy to generate at 
large values of cell sizes $\eta$ relative to $\epsilon$ and hence detect convergence rates. 

\begin{figure}[h]       
    \centering
 \includegraphics[width=0.45\textwidth]{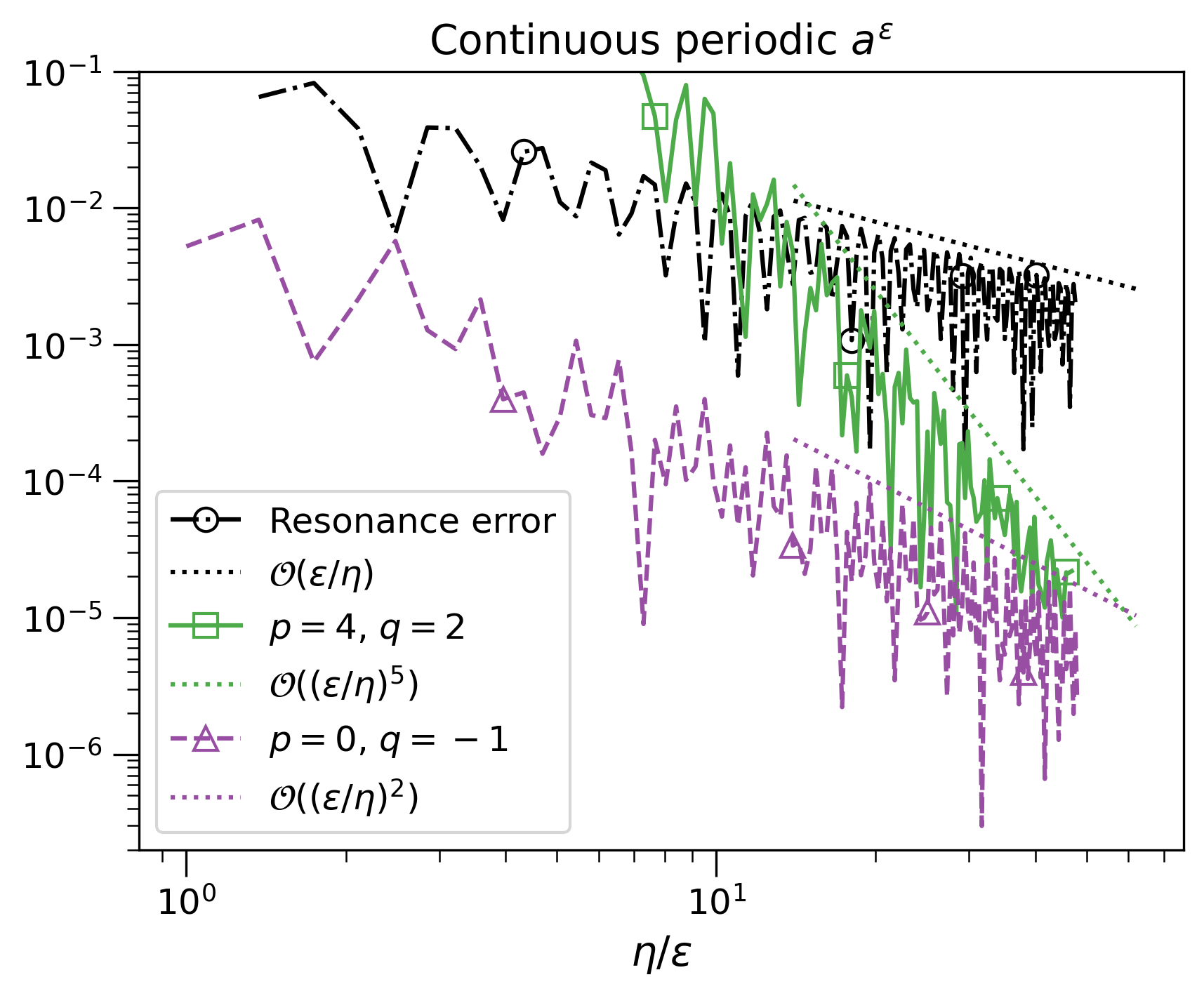}
 \includegraphics[width=0.45\textwidth]{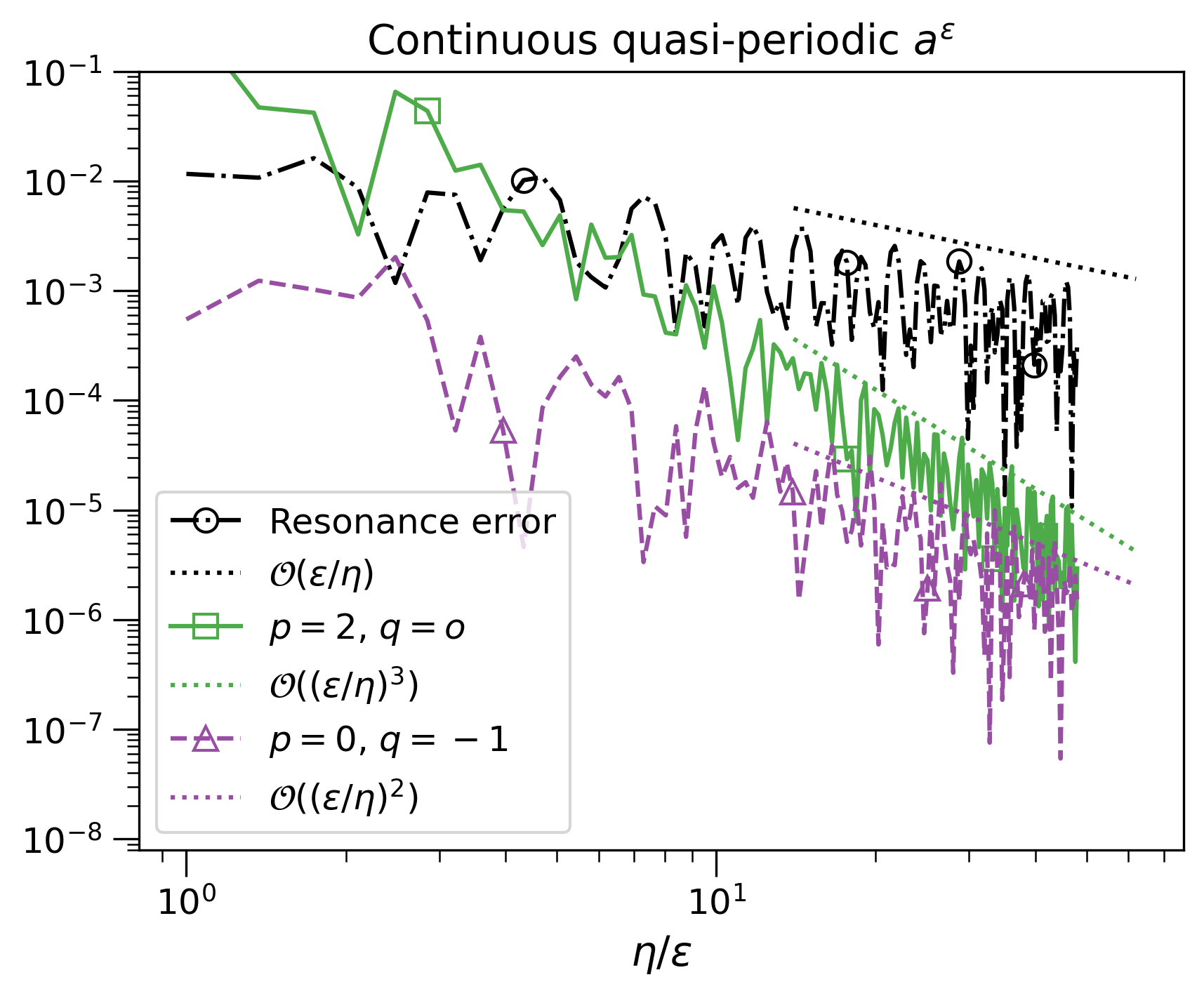}
 \includegraphics[width=0.45\textwidth]{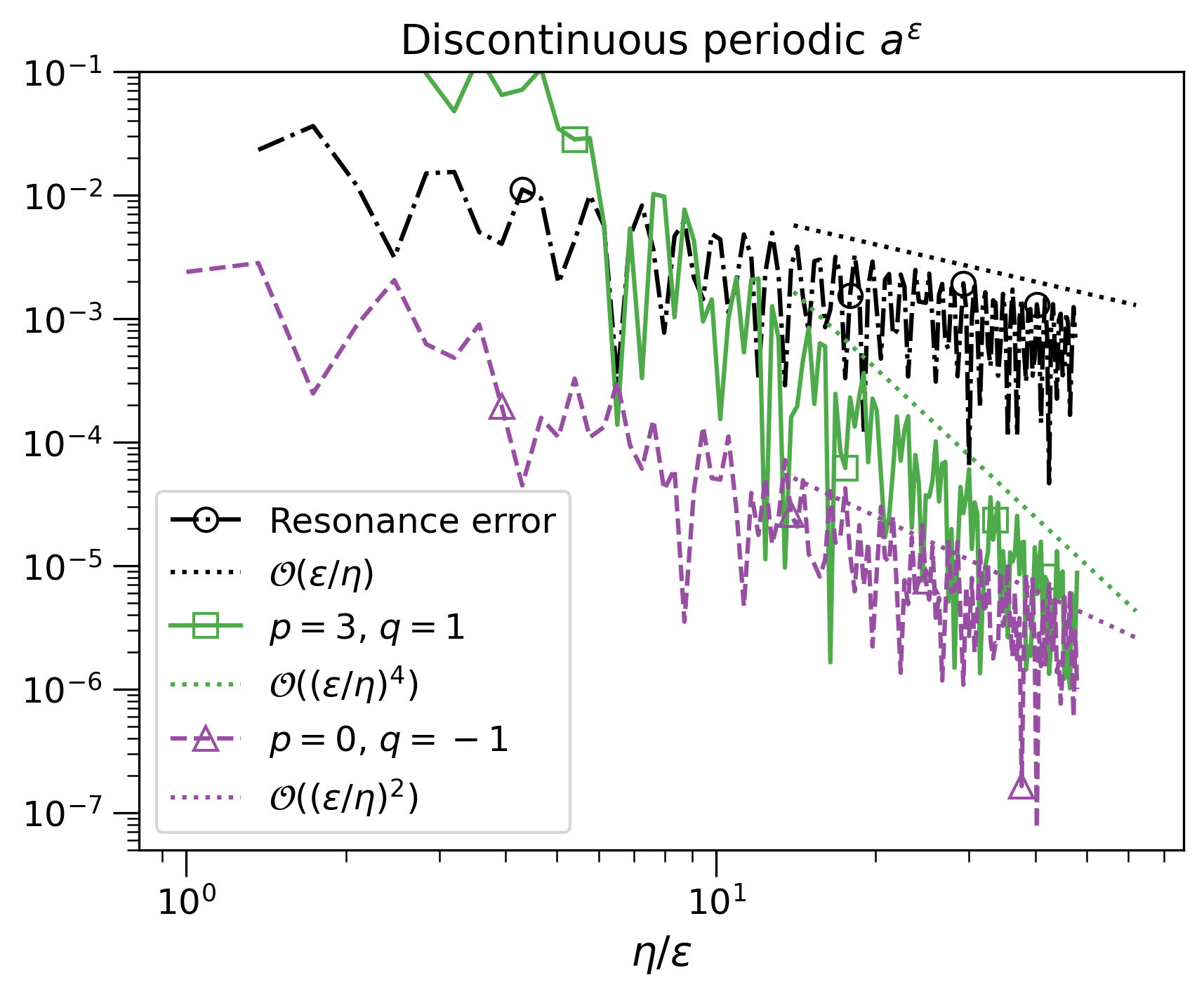}
\caption{ 
A comparison of the absolute value of the one-dimensional cell resonance error with the 
absolute value of the averaged error for various averaging kernels $K\in\mathbb{K}^{-p,q}$ 
for the three oscillatory functions $\aeps$ given in \eqref{eq:1d_aeps_formulas}. 
} 
\label{fig:1d_errors}
\end{figure}

\cref{fig:1d_errors} compares the cell resonance error 
(given by the difference between \eqref{eq:naive_avg} and the 
true harmonic average $\overline{a}$) with the averaged
error \eqref{eq:smooth_avg} for various kernels $K \in \mathbb{K}^{-p,q}$.
We consider three oscillatory functions $\aeps(x)$ given by 
\begin{align}
\aeps_1(x) =&(1.1 + \cos(2\pi x/\epsilon))^{-1} \nonumber \\
\aeps_2(x) =&(2.2 + \cos(2\pi x/\epsilon) + \cos(\sqrt{2} 2 \pi x/\epsilon))^{-1} \nonumber \\
\aeps_3(x) =& (2 + \sign (\cos(2\pi x/\epsilon)))^{-1}.  \label{eq:1d_aeps_formulas}
\end{align}
The first two examples are smooth, while the third is discontinuous. In all cases
the true harmonic average, and, more generally, the incorrect average $\rhoeps(\eta)$ are 
easily computed analytically; note that the second example is quasi-periodic, and hence the true 
harmonic average is given by \eqref{eq:naive_avg} as $\eta \to\infty$. The composite
trapezoidal rule is used to compute the weighted averages \eqref{eq:smooth_avg}. We use a large
number of quadrature points to ensure the asymptotic decay rates are computed accurately;
in particular when $\eta=\epsilon$, $N_{\rm trap}= 4096$ points are used, and the number of points
is increased proportionally as $\eta$ increases. All kernels $K \in \mathbb{K}^{-p,q}$ that are 
used are of the form \eqref{eq:kernel_form}. 

In each case, the cell resonance error of course decays as $\epsilon/\eta$, while the 
averaged errors decay at the rates predicted by \cref{thm:1d_err_thm} and 
\cref{thm:second_order}. Although the convergence rates in $\epsilon/\eta$ are higher for
averaging kernels with larger $p$ and $q$, the constant in front of the asymptotic decay 
rate also increases. Empirically we find that $\norm{K}_{\infty}$ increases with $p$ 
and $q$, which impacts the estimate \eqref{eq:kernel_decay_estimate} in  
\cref{lem:avg_kernels2}. 
In contrast, the simple arithmetic average (using $K(x)=1$) only decreases as 
$(\epsilon/\eta)^2$ as predicted by \cref{thm:second_order}. However, the 
asymptotic constant is relatively small, so it still 
results in a larger decrease of the resonance error even at cell sizes $\eta \sim 50\epsilon$; 
\emph{a fortiori}, it performs the best at ``intermediate'' values of $\eta$ 
(e.g.\ $3\epsilon \lesssim \eta \lesssim 10\epsilon$) that are more practical in higher dimensions.

\subsection{Dimension $d=2$}\label{subsec:results2d}
We next present two dimensional numerical computations of the cell resonance error
as a function of domain size $\eta$ that suggest something analogous to the 
one-dimensional decomposition from \cref{lem:expansion} indeed holds in higher
dimensions. 


The cell problems \eqref{eq:chi_j_problem} on $I_{\eta} = [-\eta/2,\eta/2]^2$
are solved with a second order finite difference method implemented in C++ with 
the linear algebra library Eigen \cite{eigenweb}.
$N= 32$ grid points are used 
in both the $x_1$ and $x_2$ directions for the case $\eta = \epsilon$; 
as $\eta$ is increased, $N$ is increased proportionally so that the resolution 
is approximately constant for each problem. The estimates for the entries in 
the homogenized tensor $\rho_{ij}(\eta)$ (given by \eqref{eq:general_rho_ij}) 
are computed with the composite trapezoidal rule, where the first order partial 
derivatives in the integrand are approximated with a second order centered
 difference formula. The estimates are computed for a sequence of domain sizes 
between $\eta=\epsilon$ and $\eta = 16\epsilon$ with a spacing of 
$\Delta \eta = 0.05\epsilon $. 

\begin{figure}[h]       
    \centering
 \includegraphics[width=0.45\textwidth]{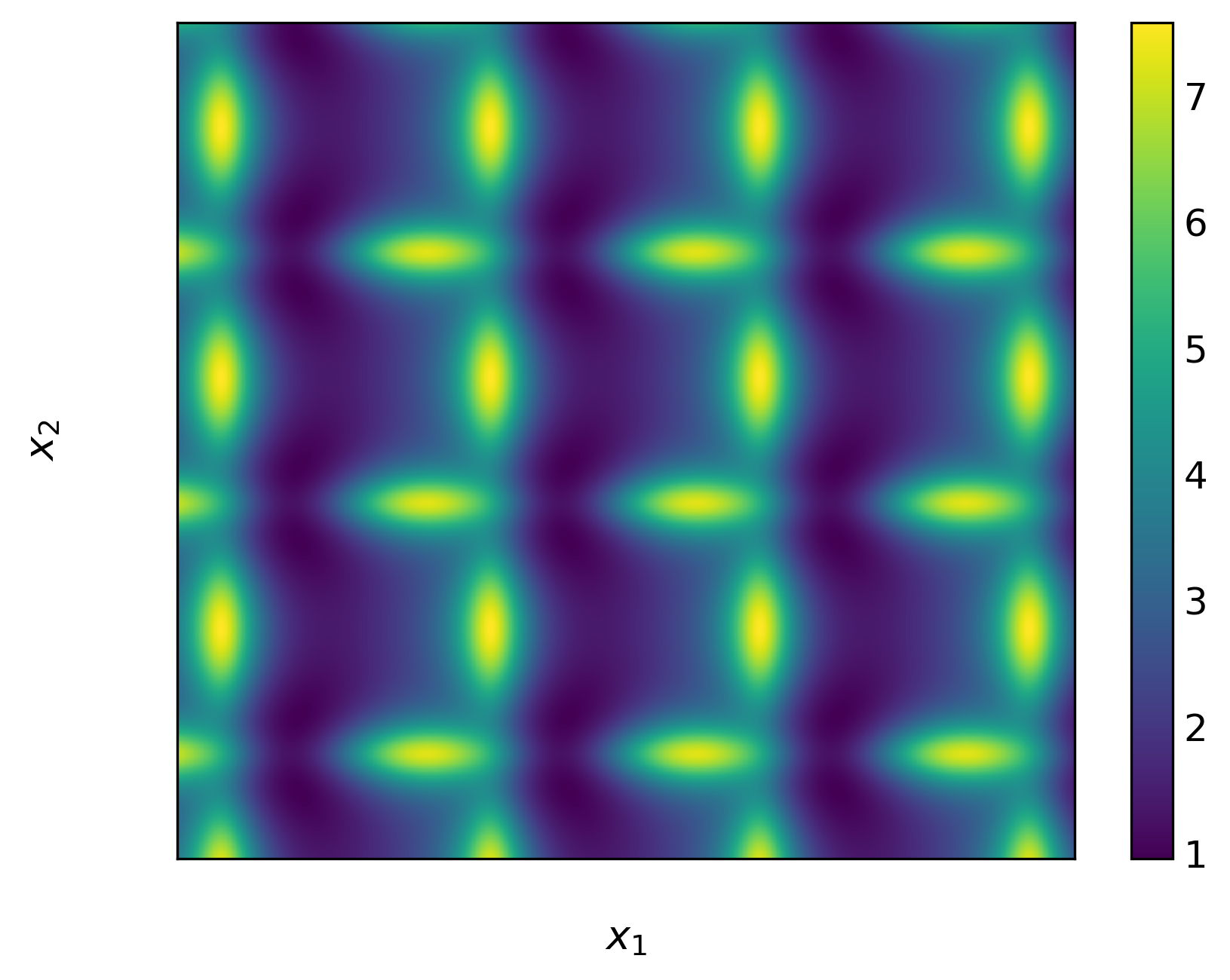}
 \includegraphics[width=0.45\textwidth]{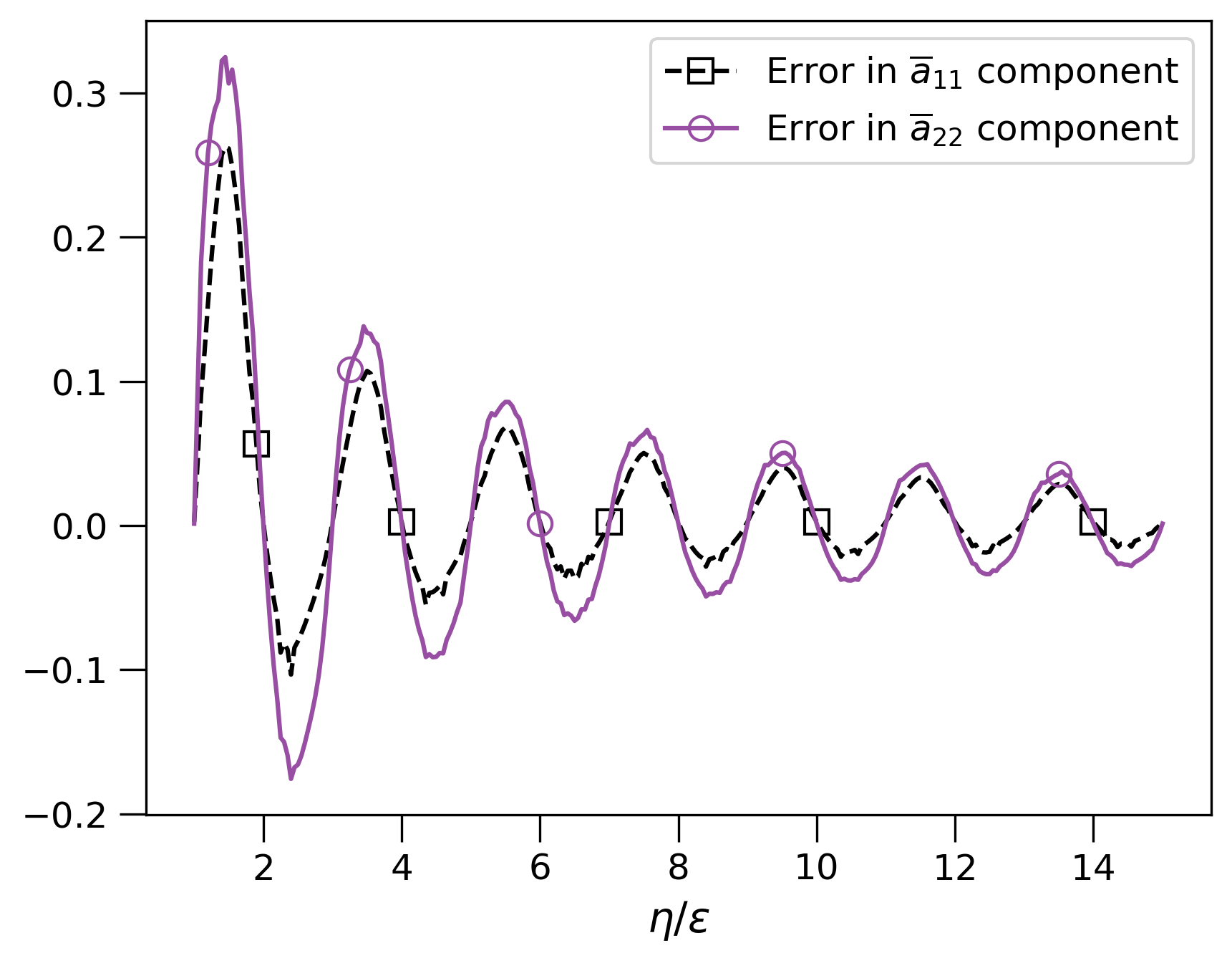}
\caption{ 
Visualization of the periodic oscillatory function \eqref{eq:aeps_case2} (left) and 
corresponding components of the cell resonance error as a function of domain 
size (right).  
} 
\label{fig:tensor_visuals_and_resonance_errors2}
\end{figure}

We consider isotropic oscillatory tensors; that is, $a^{\epsilon}$ is given by 
a (strictly positive) scalar-valued function times the identity matrix. 
Consider first the $\epsilon$-periodic function 
\begin{equation}  \label{eq:aeps_case2}
a^{\epsilon}(x_1,x_2) =  \frac{2+1.5 \sin(2\pi x_1/\epsilon)}{2 + 1.5 \sin(2\pi x_2/\epsilon)} 
+ \frac{2+1.5 \sin(2\pi x_2/\epsilon)}{2 + 1.5 \sin(2\pi x_1/\epsilon)},
\end{equation} 
taken from Section 3.3~in \cite{yue:2006} and visualized in  
\cref{fig:tensor_visuals_and_resonance_errors2} (left). In this case the true homogenized 
tensor $\overline{a}$ is diagonal but anisotropic. The entries are approximated by 
solving the cell problem \eqref{eq:chi_j_problem} for $\eta=\epsilon$ with $N=1024$ 
grid points in both the $x_1$ and $x_2$ directions, which gives
$$
\overline{a} \simeq
\begin{pmatrix}
 2.34348520086 & 0  \\ 
0 & 2.87329450077
\end{pmatrix}.
$$
The resonance error $\rho_{ij}(\eta) - \overline{a}_{ij}$ for the two entries 
$i=j=1$ and $i=j=2$ is also shown in \cref{fig:tensor_visuals_and_resonance_errors2} (right); 
between the $\epsilon/\eta$ decay envelopes the error is clearly oscillatory. 

\begin{figure}[h]       
    \centering
 \includegraphics[width=0.45\textwidth]{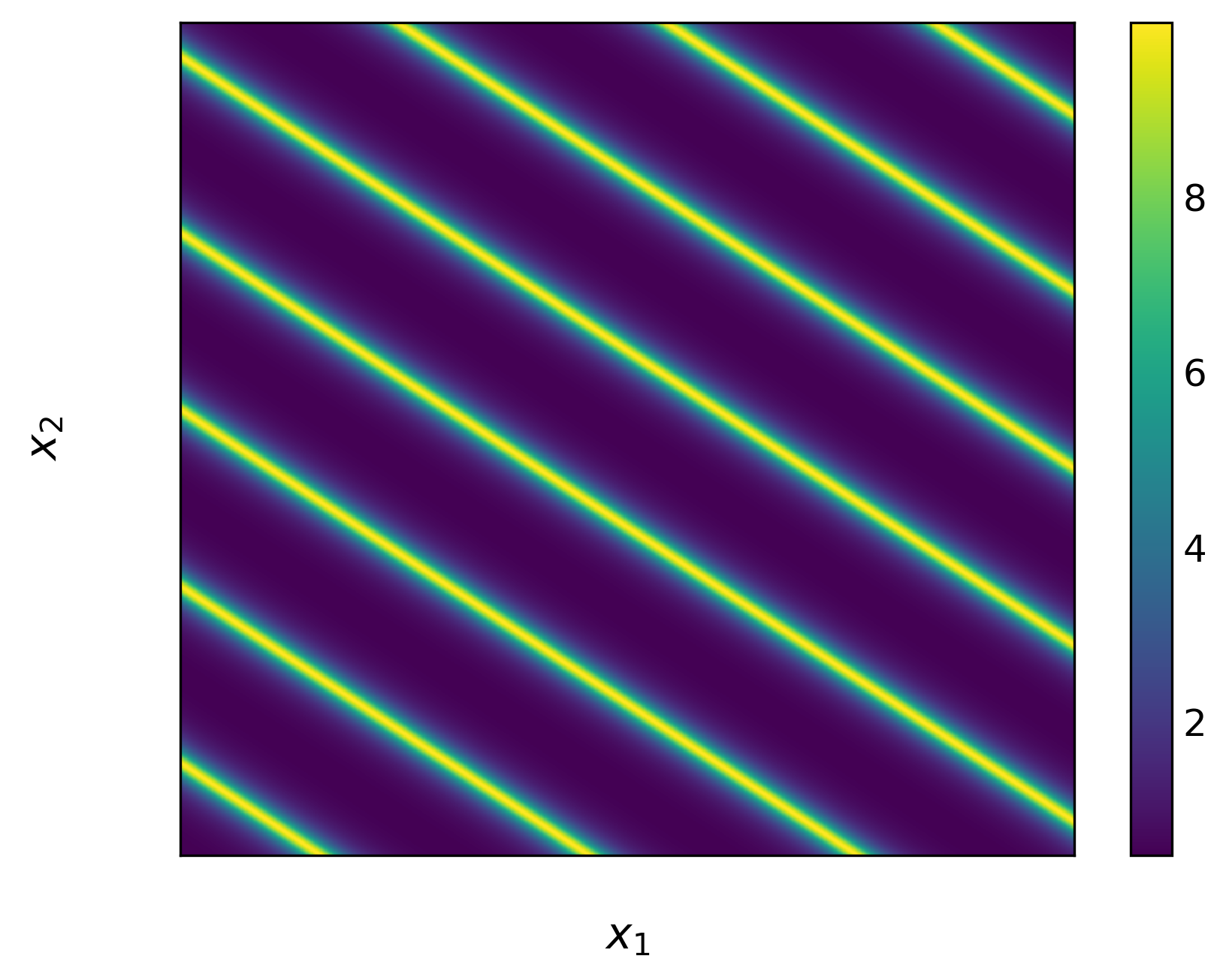}
 \includegraphics[width=0.45\textwidth]{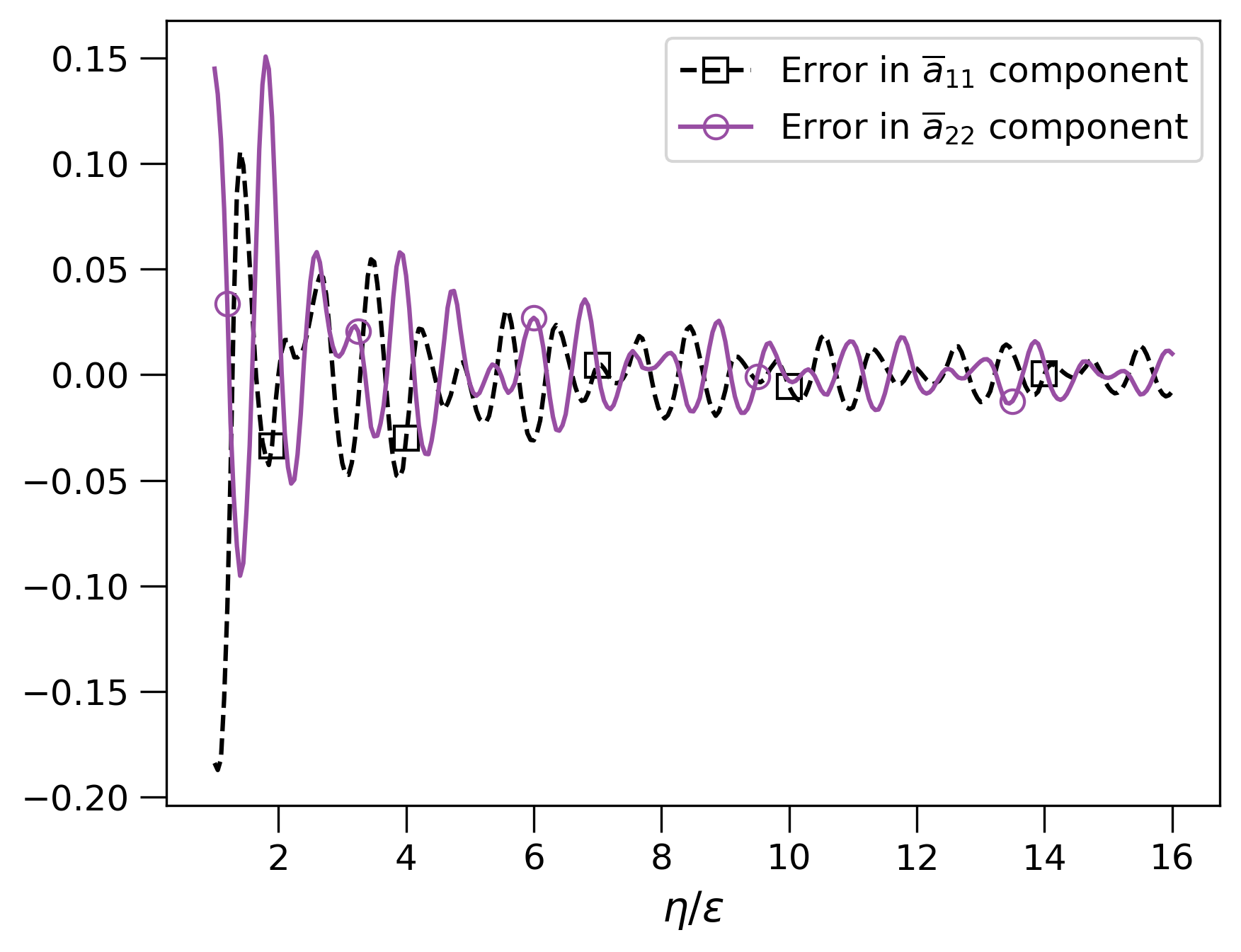}
\caption{ 
Visualization of the quasi-periodic oscillatory function \eqref{eq:aeps_case4} (left) and 
corresponding components of the cell resonance error as a function of domain 
size (right).  
} 
\label{fig:tensor_visuals_and_resonance_errors4}
\end{figure}

\cref{fig:tensor_visuals_and_resonance_errors4} (right) shows that the 
resonance error is also oscillatory for the isotropic tensor given by the identity 
matrix times the quasi-periodic function  
\begin{equation}  \label{eq:aeps_case4}
\aeps(x_1,x_2) = \big(1.1 + \cos(2\pi(x_1 + \sqrt{2} x_2)/\epsilon)\big)^{-1}. 
\end{equation} 
As in the the first case, the true homogenized tensor is anisotropic; in this
case however the coefficients $\overline{a}_{11}$ and $\overline{a}_{22}$ are 
approximated by solving \eqref{eq:chi_j_problem} at relatively large values of $\eta/\epsilon$. 
In particular, we take coefficients to be the arithmetic average of the values 
$\rho_{ij}(\eta)$ from $\eta =13\epsilon$ to $\eta = 16\epsilon$ which gives
$$
\overline{a} \simeq
\begin{pmatrix}
 1.75643523765 & 0  \\ 
0 & 1.34396605902 
\end{pmatrix}.
$$

Given the oscillatory nature of the cell resonance error as a function of domain 
size $\eta$, we now apply the numerical strategy outlined in 
\cref{subsec:algo_description} for its reduction to the cases \eqref{eq:aeps_case2} 
and \eqref{eq:aeps_case4}. In both cases, the weighted average \eqref{eq:S_ij} is 
approximated by the composite trapezoidal rule with intervals of size $\Delta t = 0.05 \epsilon$. 

\begin{figure}[h]       
    \centering
 \includegraphics[width=0.45\textwidth]{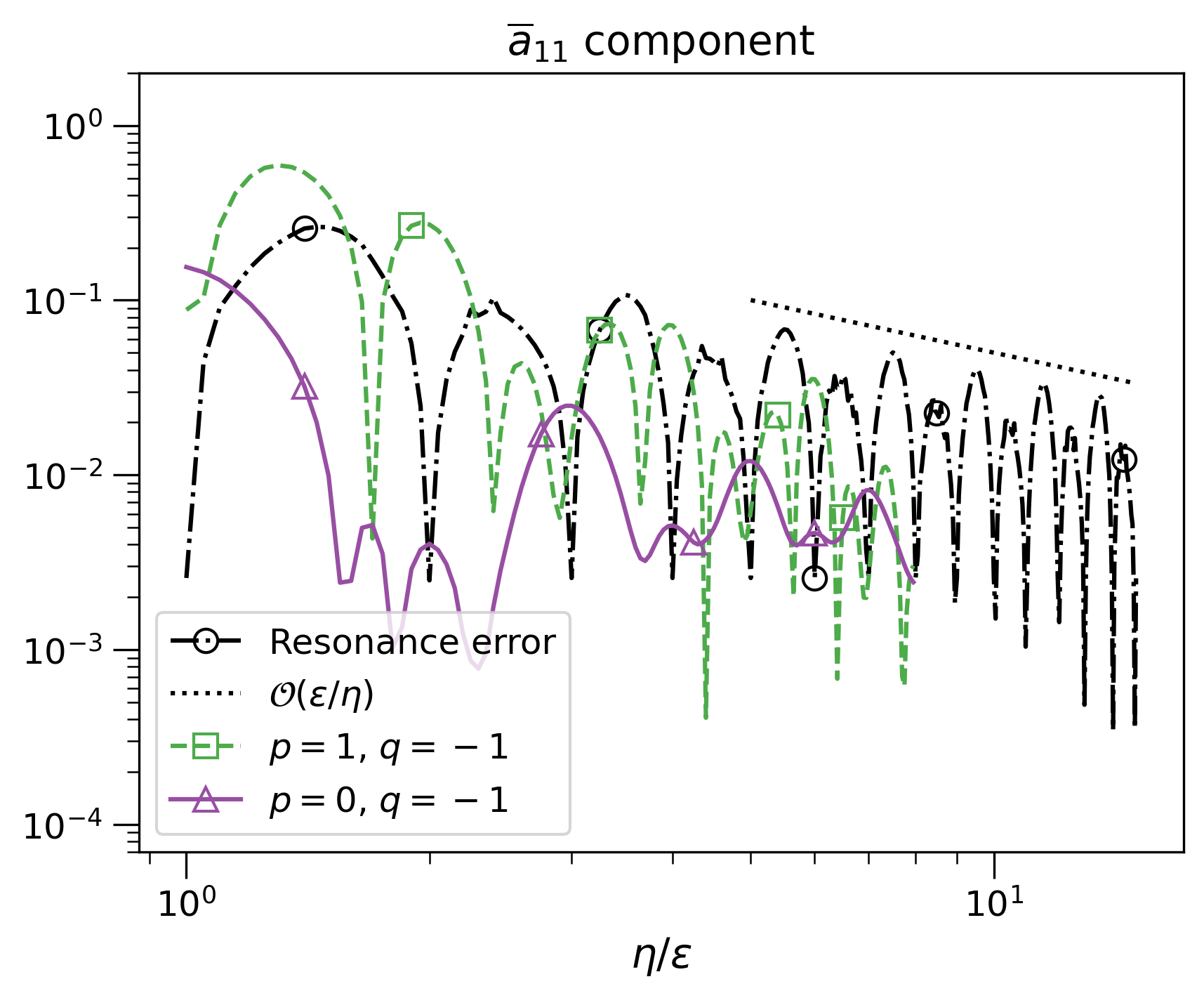}
 \includegraphics[width=0.45\textwidth]{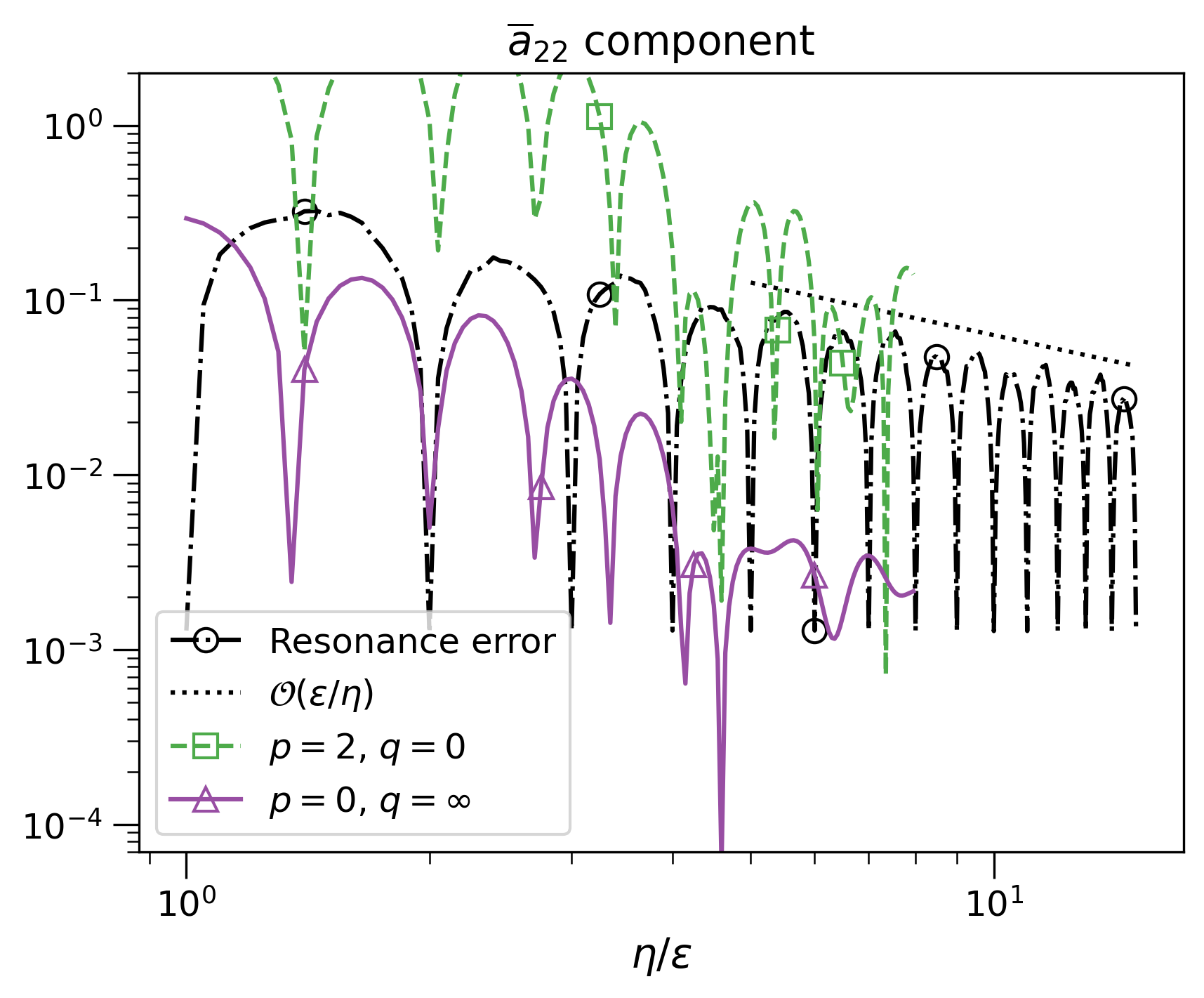}
\caption{ 
A comparison of the absolute value of the resonance error with the 
absolute value of the averaged error for various averaging kernels 
$K\in\mathbb{K}^{-p,q}$ for the two-dimensional case \eqref{eq:aeps_case2}.
} 
\label{fig:smoothed_err_2}
\end{figure}

\cref{fig:smoothed_err_2} shows the results for each of the two components 
of the resonance error when $\aeps$ is given by the $\epsilon$-periodic function \eqref{eq:aeps_case2}.
We consider several types of kernels $K \in \mathbb{K}^{-p,q}$. As discussed
in \cref{sec:numerical_results1d} for the one dimensional results, 
kernels with $p \ge 1$ lead to unsatisfactory results at the relatively low values
of $\eta/\epsilon$ relevant for dimensions $d\ge 2$ owing to their large $\infty$-norms. In 
contrast, the kernels $K(x) = 1 \in \mathbb{K}^{0,-1}([1,2])$ and 
\begin{equation}\label{eq:K0exp}
 K(x) = C e^{1/((x-1)(x-2))} \in \mathbb{K}^{0,\infty}([1,2])
\end{equation}
(where $C$ is chosen so that $K$ has unit mass on $[1,2]$) reduce the resonance
error by an order of magnitude or more at relatively low values of $\eta/\epsilon$.

\begin{figure}[h]       
    \centering
 \includegraphics[width=0.45\textwidth]{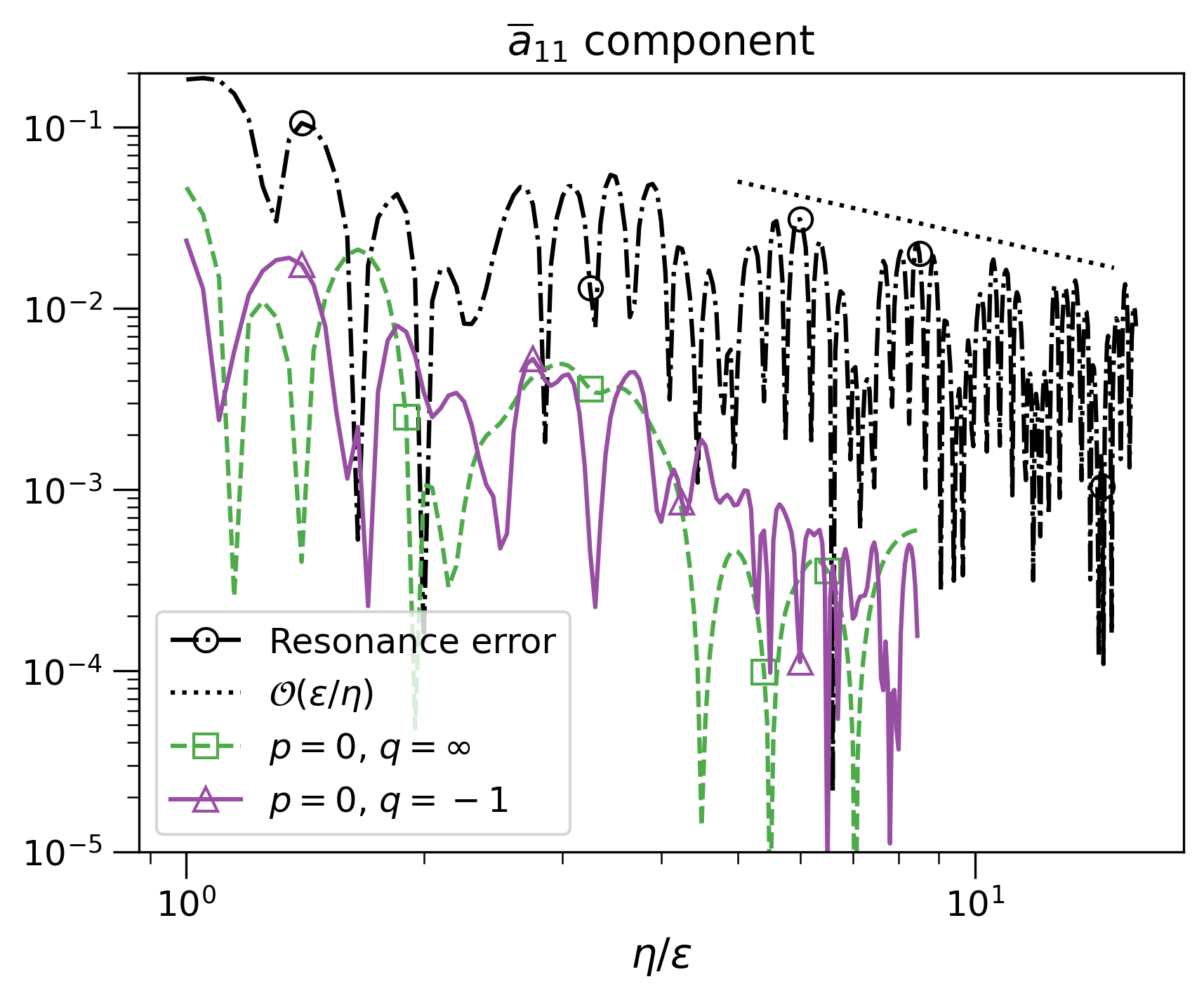}
 \includegraphics[width=0.45\textwidth]{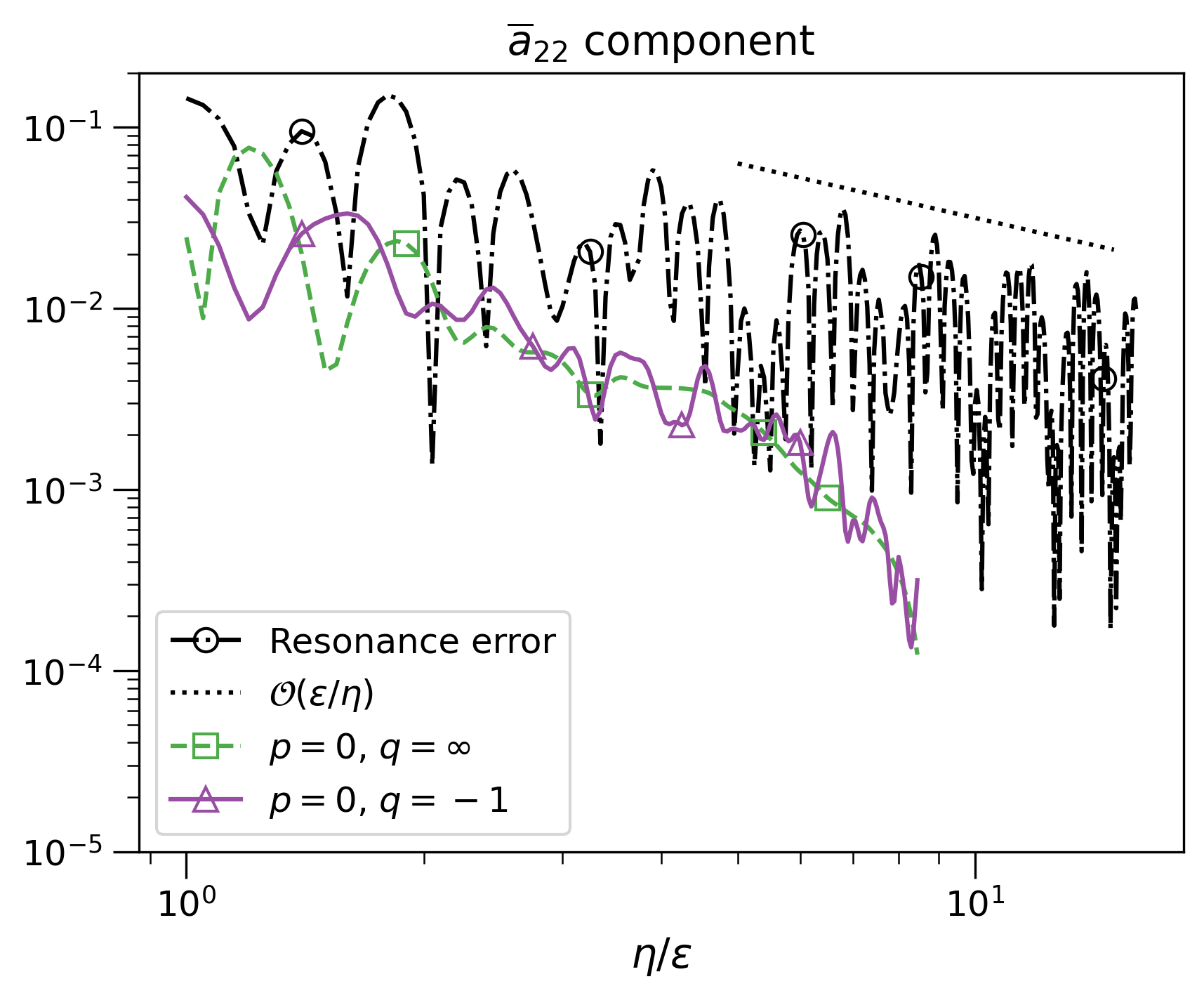}
\caption{ 
A comparison of the absolute value of the resonance error with the 
absolute value of the averaged error for various averaging kernels 
$K\in\mathbb{K}^{-p,q}$ for the two-dimensional case \eqref{eq:aeps_case4}.
} 
\label{fig:smoothed_err_4}
\end{figure}

\cref{fig:smoothed_err_4} shows the results for each of the two components 
of the resonance error when $\aeps$ is given by the quasi-periodic function \eqref{eq:aeps_case4}. In this 
case we only consider the kernels $K(x)=1$ and the exponential kernel \eqref{eq:K0exp}. 
As in the previous example, both reduce the resonance error by several 
orders of magnitude.

\subsection{Computational cost and discussion}
\label{sec:discussion}
Assuming the estimate for the cell resonance error in \cref{thm:1d_err_thm} generalizes to 
higher dimensions,
the computational cost of the method proposed in  \cref{subsec:algo_description}  
to reduce the error below some tolerance is asymptotically
similar to the cost of the hyperbolic microsolver proposed in \cite{runborg:2016}. 
Indeed, assume a linear system with $N$ degrees
of freedom can be solved with computational cost proportional to $N^{\alpha}$; for 
example, $\alpha = 2$ (respectively $2.5$) for a sparse LU solver in dimension $d=2$
(resp.\ $d=3$). In an ideal scenario, 
$\alpha=1$ could be achieved with a multigrid method, however, for microscale domains with
diameter $\eta$ sufficiently large relative to $\epsilon$, it is not clear that the 
coarse grid solvers can resolve the oscillatory nature of the solution.  
To generate the approximation \eqref{eq:S_ij} to the true homogenized coefficients, 
the method proposed here consists of solving a sequence of elliptic problems on domains
with diameters that range from $\eta$ to $2\eta$. Since \eqref{eq:general_rho_ij} can 
oscillate with period $\epsilon$, in general $\bigoh{\eta/\epsilon}$ solves are required. 
Similarly, for the hyperbolic solver a time step $\Delta t \sim \epsilon$ is required
for a total time of integration $T \sim \eta$. Hence, assuming that the total computational 
cost to estimate the homogenized coefficients (including, e.g.\ the use of quadrature formulas) 
is dominated by linear solves, both methods have a total cost proportional to 
$(\eta/\epsilon)^{\alpha d + 1}$. 
Of course, to be fair it should again be emphasized that the authors 
in \cite{runborg:2016} rigorously showed that the resonance error can be made to decay 
like $(\epsilon/\eta)^s$, where $s$ can be arbitrarily large, while the same can only be 
said for the current approach in dimension $d=1$.

Of all the existing techniques in the literature to reduce the cell resonance error, 
in theory 
the modified elliptic problem recently proposed in \cite{abdulle2023elliptic}  
ought to be the most computationally efficient,  
assuming the Arnoldi iterations used to generate the forcing function $f$ 
(see Eq.~\eqref{eq:abdulle_forcing}) converge sufficiently rapidly; see section 4.3 in
\cite{abdulle2023elliptic} for a discussion of this point. As with other analyses 
found in the literature (see, for example, section 4.1 of \cite{runborg:2016}), however, 
this 
assumes that
the cell domain size $\eta$ is large enough that the theoretical decay rates saturate, which 
may not be true until $\eta \sim 10\epsilon$ or larger. 
In practice a microsolver's utility is determined by its performance at 
relatively low values of $\eta$; depending on the application, this could mean $\eta$ is several 
times the assumed local period, or correlation length, $\epsilon$ of the heterogeneous media.
A comparison of the existing 
microsolvers
in the literature at such ``practical" values of $\eta$ is beyond the scope
of the present work, however such a study may be valuable for streamlining research in 
numerical homogenization, 
especially if done in collaboration with the engineering community \cite{bargmann2018generation}. 

One advantage of the present approach
is that the microscale domain calculations 
do not communicate with one another, i.e.\ the microscale solution at a given domain size 
does not depend on any of the others. The efficiency of the method hence can be 
greatly increased if the solutions are computed independently in parallel. 

Another feature unique to the present approach is that the form of the microscale problem 
is unmodified from the standard elliptic problem from homogenization theory; there is hence
no need to develop approximations for other operators. 
It is natural then to combine the approach with reduced basis (RB) techniques for numerical 
homogenization, see for example 
\cite{abdulle2014reduced,abdulle2015reduced,boyaval2008reduced,nguyen2008multiscale}.
%
Such techniques are important for three dimensional problems, for which the cost of solving 
microscale problems 
at every macroscale grid point where homogenized coefficients are needed
can be prohibitively expensive. They are also important for optimal control problems constrained 
by multiscale equations of the form \eqref{eq:ms_elliptic}. Resolving the boundary resonance error in this context is an open problem, 
as noted in \cite{abdulle2012reduced}. 

Broadly speaking, RB methods precompute in an ``offline'' stage a low dimensional set of basis functions
(say of dimension $N$) by solving
\eqref{eq:chi_j_problem} at different locations throughout the computational domain $\Omega$. 
Homogenized coefficients needed by a macrosolver are then efficiently computed in an ``online'' stage by 
solving linear systems of size $N\times N$ in the RB space. 
The strategy described in \cref{subsec:algo_description} could be used to reduce the resonance error for RB methods by 
first precomputing $M$ different basis sets, each corresponding to a different domain size $t_m$, 
$1 \le m \le M$.  This of course incurs a larger computational 
cost in the offline stage, although we again note that this process is parallelizable. 
Once the different bases are computed, 
however, the values \{$\rho_{ij}(t_m)\}_{m=1}^M$ in \eqref{eq:final_quad} can be easily obtained by 
solving $M$ low dimensional linear systems. We plan to explore this possibility in future research. 

\section{Conclusions}
\label{sec:conclusions}
A novel method for reducing the boundary resonance error in numerical homogenization is proposed.
Rather than modifying the form of the microscale problem itself, the 
method is based on taking a weighted average of the boundary error at a 
sequence of different microscale domain sizes.  

Underlying the method is the observation that
the boundary error itself is an oscillatory
function of domain size.
%
%
In one dimensional and two dimensional
``tubular'' domains the oscillatory nature is rigorously characterized. Numerical 
evidence suggests the results also hold in more general settings, however,
more work is required to rigorously show that this is indeed the case. The problem 
can be cast as a classical elliptic homogenization problem with oscillatory 
boundary data, as considered in \cite{masmoudi:2012} and \cite{feldman:2014}; 
understanding the potentially complicated boundary layers discussed 
therein is key for further progress in this direction.

Based on the oscillatory nature of the boundary error, we propose 
solving microscale problems for a sequence of domain sizes and then  
averaging the results against kernels 
with vanishing ``negative'' moments designed to 
accelerate the convergence of the resonance error to zero at infinitely large 
microscale domains. 

Owing to the rigorous analysis of the behavior of the resonance error, we show for dimension 
$d=1$ that the error can be made to decay as 
$(\epsilon/\eta)^s$ for $s$ arbitrarily large, which we confirm with numerical 
experiments. However, the method is generally designed for higher dimensions, and numerical 
examples in $d=2$ demonstrate its efficacy. 

Because the form of the microscale problem is unchanged from classical 
homogenization theory, the method does not require developing approximations for other operators. 
In future work, the method could be combined with reduced
basis numerical homogenization techniques used for three dimensional problems and 
for optimization problems with multiscale partial differential equations as constraints. 

\section*{Acknowledgments}
The authors acknowledge support from the Oden Institute for Computational Engineering
and Sciences. We also thank Eric Hester for assistance with Mathematica calculations
to determine the averaging kernels used in some of the numerical examples.

\appendix
\section{Proof of results for averaging kernels} \label{sec:appendix}
\textit{Proof of \cref{lemma:primitives}:}
By continuity of $g$, both the left and right-handed 
limits exist for all $x$: 
$$\left| g^{[1]}(x\pm \Delta x) - g^{[1]}(x) \right| \to 0$$ 
as $\Delta x \to 0$, so that $g^{[1]}$ is continuous. Furthermore, 
\begin{align*}
g^{[1]}(x+1) &= \int_0^{x+1} g(s)\, ds + c_1 = \int_0^x g(s) \, ds + c_1 + \int_x^{x+1} g(s) \, ds  \\
&= g^{[1]}(x) + \int_0^1 g(s) \, ds = g^{[1]}(x)
\end{align*}
from the assumption that $g$ is mean zero. Lastly, note that continuous, 
periodic functions are bounded. $\square$


The next proof is an adaptation of results from \cite{runborg:2014}.

\textit{Proof of \cref{lem:avg_kernels2}:}
Define $F^{[0]}(x) = f(x) - \chevron{f}$, and let $F^{[q+1]}$ be the $q+1$ primitive, 
as in \cref{dfn:primitive}. Note 
\begin{equation*}
\frac{d^{q+1}}{dx^{q+1}} F^{[q+1]}(x) = F^{[0]}(x),
\end{equation*}
so that 
\begin{equation}\label{eq:iterated_derivs}
F^{[0]}(x\eta/\epsilon) = (\epsilon/\eta)^{q+1} \frac{d^{q+1}}{dx^{q+1}} F^{[q+1]}(x\eta/\epsilon).
\end{equation}
Furthermore, note that $\norm{F^{[q+1]}}_{\infty} \le 2^{q+1} \norm{f}_{\infty}$. Then after a change of variables,
\begin{equation*}
\int_{\eta}^{2\eta} K_{\eta}(x) x^{-r} f\left(\frac{x}{\epsilon}\right) dx = 
\underbrace{\eta^{-r} \int_{1}^{2} K(u) u^{-r} F^{[0]}(u\eta /\epsilon) du}_{(\textrm{I})} + 
\underbrace{\chevron{f} \eta^{-r} \int_{1}^{2} K(u) u^{-r}  du,}_{(\textrm{II})} 
\end{equation*}
where $F^{[0]}$ is mean-zero. Since $K \in \mathbb{K}^{-p,q}$, $(\textrm{II})$ vanishes
if $r \le p$; otherwise, 
\begin{equation*}
\left|\chevron{f} \eta^{-r} \int_{1}^{2} K(u) u^{-r} du \right| \le\left|\chevron{f}\right| \norm{K}_{\infty} \, \eta^{-r}.
\end{equation*}
After using \eqref{eq:iterated_derivs} and integrating by parts, $(\textrm{I})$ becomes
\begin{align*}
(\textrm{I}) &= (\epsilon/\eta)^{q+1} \eta^{-r} \int_{1}^{2} K(u) u^{-r}\frac{d^{q+1}}{du^{q+1}} F^{[q+1]}(u\eta/\epsilon) du \\  
&= (-1)^{q+1} (\epsilon/\eta)^{q+1} \eta^{-r}\int_{1}^{2} \frac{d^{q+1}}{du^{q+1}}\left( K(u) u^{-r}\right) F^{[q+1]}(u\eta/\epsilon) du 
\end{align*}
Since $u^{-r} \in C^{\infty}([1,2])$ and by assumption $K^{(q+1)} \in BV([1,2])$, 
$$ g(u):= \frac{d^{q+1}}{du^{q+1}}\left( K(u) u^{-r}\right)  \in BV([1,2])$$
as well. Changing variables again and breaking the resulting integral in three pieces gives
\begin{align*}
(\textrm{I}) &=   
(-1)^{q+1} (\epsilon/\eta)^{q+2} \eta^{-r}\int_{\eta/\epsilon}^{2\eta/\epsilon} g(t\epsilon/\eta) F^{[q+1]}(t) \, dt  \\
&= (-1)^{q+1} (\epsilon/\eta)^{q+2} \eta^{-r} 
\Big( \int_{\eta/\epsilon}^{\ceil{\eta/\epsilon}} g(t\epsilon/\eta) F^{[q+1]}(t) \, dt
+ \int_{\ceil{\eta/\epsilon}}^{\floor{2\eta/\epsilon}} g(t\epsilon/\eta) F^{[q+1]}(t) \, dt \\
&\,\,\,\,\qquad\qquad\qquad\qquad\qquad +\int_{\floor{2\eta/\epsilon}}^{2\eta/\epsilon} g(t\epsilon/\eta) F^{[q+1]}(t) \, dt \Big)
\end{align*}
The first integral in parentheses can be bounded as 
$$ 
\Big| \int_{\eta/\epsilon}^{\ceil{\eta/\epsilon}} g(t\epsilon/\eta) F^{[q+1]}(t) \, dt \Big| 
\le \norm{g}_{\infty} \norm{F^{[q+1]}}_{\infty} \le 2^{q+1} \norm{f}_{\infty} \,  V_1^2(g)  
$$
where $V_1^2(g)$ denotes the total variation of $g$ on $[1,2]$ and depends only on $K$, its derivatives, 
and $q$. A similar bound exists for the 
third integral. Next, let $N = \floor{2\eta/\epsilon} - \ceil{\eta/\epsilon}$. By construction, $F^{[q+1]}$ 
is a mean-zero, 1-periodic function, so that 
$$
\int_{\ceil{\eta/\epsilon}}^{\floor{2\eta/\epsilon}}  F^{[q+1]}(t) \, dt =0 .
$$
The second integral is then bounded by 
$$
\Big| \sum_{j=1}^{N-1} \int_{\ceil{\eta/\epsilon}j}^{\ceil{\eta/\epsilon}(j+1)} 
\big( g(t\epsilon/\eta) - g(\ceil{\eta/\epsilon}j\epsilon/\eta) \big) F^{[q+1]}(t) \, dt\Big|  
\le 2^{q+1} \norm{f}_{\infty} \,  V_1^2(g).
$$
Putting it all together, we've shown that
$$
\Big| \int_{\eta}^{2\eta} K_{\eta}(x) x^{-r} f\left(\frac{x}{\epsilon}\right) dx \Big| 
\le 3 \eta^{-r} \left(\frac{\epsilon}{\eta}\right)^{q+2}  2^{q+1} \norm{f}_{\infty} V_{1}^2(g) + 
\left|\chevron{f}\right| \norm{K}_{\infty} \, \eta^{-r}
$$
for $p < r$ and 
$$
\Big| \int_{\eta}^{2\eta} K_{\eta}(x) x^{-r} f\left(\frac{x}{\epsilon}\right) dx \Big| 
\le 3 \eta^{-r} \left(\frac{\epsilon}{\eta}\right)^{q+2}  2^{q+1} \norm{f}_{\infty} V_{1}^2(g)  
$$
for $p \ge r$, as desired. $\square$


\section{Proof of \cref{thm:second_order}} \label{sec:appendix_second_order_decay}
Following the proof of \cref{thm:1d_err_thm}, we use the expansion 
\eqref{eq:main_decomposition} from \cref{lem:expansion} and then take 
the arithmetic average on the interval $[\eta,2\eta]$. Since the expansion 
converges absolutely for $\eta/\epsilon > \overline{a} B$, this means
\begin{equation}\label{eq:B1}
\frac{1}{\eta} \int_{\eta}^{2\eta} \rhoeps(x) \, dx - \overline{a} 
= \overline{a} \sum_{k=1}^{\infty} \frac{1}{\eta} \int_{\eta}^{2\eta} \left(-\frac{\epsilon}{x} q(x/\epsilon)\right)^{k} dx
\end{equation}
where 
$$ q(x/\epsilon) = \overline{a} \int_0^{x/\epsilon} (b(s)-\chevron{b})\, ds.$$
Considering first the $k=1$ term in the sum and integrating by parts gives 
\begin{align}
-\frac{1}{\eta} \int_{\eta}^{2\eta} \frac{\epsilon}{x} q(x/\epsilon) \, dx 
&= -\frac{1}{\delta} \int_{\delta}^{2\delta} \frac{q(u)}{u}  \, du \nonumber \\
&= \frac{1}{\delta} \int_{\delta}^{2\delta} \frac{Q(u)}{u^2}  \, du 
- \frac{1}{\delta} \left( \frac{1}{2\delta} Q(2\delta) - \frac{1}{\delta} Q(\delta) \right),  \label{eq:IBP}
\end{align}
where $\delta := \eta/\epsilon$ and $Q$ is the primitive of $q$ as defined by \cref{dfn:primitive}. 
Since $b(u)-\chevron{b}$ is
a mean-zero periodic function, \cref{lemma:primitives} guarantees that $q(u)$ is also 
a periodic function. Moreover, since $b(s)=1/a(s)$ is assumed to be an even function, $q(u)$ 
is an odd function, and hence
\begin{equation*}
\int_{-1/2}^{1/2} q(u) \, du = 0. 
\end{equation*}
Because $q$ is a mean-zero, one-periodic function, \cref{lemma:primitives} implies that $Q$
is also periodic, and hence bounded independently of $\delta$. Ergo, by the triangle inequality 
\eqref{eq:IBP} implies 
\begin{equation}\label{eq:leading_term}
\left| \frac{1}{\eta} \int_{\eta}^{2\eta}  \frac{\epsilon}{x} q(x/\epsilon) \, dx  \right| 
\le  2 \norm{Q}_{\infty} \left( \frac{\epsilon}{\eta}\right)^2 
\end{equation}
It then remains to consider the $k \ge 2$ terms of the sum in \eqref{eq:B1}. 
As in \eqref{eq:bound_antideriv_of_b} and \eqref{eq:Gr_bound}, 
$$
\sup_{x \in \mathbb{R} } \left| q(x/\epsilon)\right| \le \overline{a} B 
$$
and hence 
$$
\sup_{x \in \mathbb{R}} \left| (q(x/\epsilon) )^k \right| \le \overline{a}^k B^k .
$$
%
We then have
\begin{equation*}\label{eq:tail}
\overline{a} \left| \sum_{k=2}^{\infty} \frac{1}{\eta} \int_{\eta}^{2\eta} \left(-\frac{\epsilon}{x} q(x/\epsilon)\right)^{k} dx \right| 
\le \overline{a} \sum_{k=2}^{\infty} (\overline{a} B)^k \left( \frac{\epsilon}{\eta} \right)^k  
= \overline{a}\, (\overline{a} B)^2 \left( \frac{\epsilon}{\eta} \right)^2 \varphi\big(\overline{a}B \frac{\epsilon}{\eta}\big).
\end{equation*}
Along with \eqref{eq:leading_term}, this gives the desired result. $\square$

\bibliographystyle{siamplain}
\bibliography{references}
\end{document}